\documentclass[12pt,thmsa]{article}
\usepackage{amsfonts, amsmath}

\newtheorem{theorem}{Theorem}
\newtheorem{corollary}{Corollary}
\newtheorem{proposition}{Proposition}
\newtheorem{lemma}{Lemma}
\newtheorem{defi}{Definition}

\newcommand{\p}{\Bbb{P}}

\newcommand{\e}{\Bbb{E}}

\newcommand{\ind}{\mbox{\rm 1\hspace{-0.04in}I}}

\newcommand{\R}{\mbox{\rm I\hspace{-0.02in}R}}

\newcommand{\ed}{\stackrel{(d)}{=}}
\newcommand{\ud}{\mathrm{d}}
\newcommand{\eqdef}{\stackrel{\mbox{\tiny$($def$)$}}{=}}

\def\QED{\hfill\vrule height 1.5ex width 1.4ex depth -.1ex \vskip20pt}

\begin{document}
\hspace*{-0.5in} {\footnotesize This version March 3, 2007.}
\vspace*{0.9in}
\begin{center}
{\LARGE The upper envelope of positive self-similar Markov
processes. \vspace*{0.4in}}\\
{\large J.C. Pardo\footnote{Research supported by a grant from
CONACYT (Mexico).}\vspace*{0.2in}}
\end{center}
\noindent $^{1}$  {\footnotesize Laboratoire de Probabilit\'es et
Mod\`eles Al\'eatoires, Universit\'e Pierre et Marie Curie, 4, Place
Jussieu - 75252 {\sc Paris Cedex 05.} E-mail:
pardomil@ccr.jussieu.fr}\\

\noindent {\it Abstract:} {\footnotesize We establish integral tests
and laws of the iterated logarithm at $0$ and at $+\infty$, for the
upper envelope of positive self-similar Markov processes. Our
arguments are based on the Lamperti representation, time reversal
arguments and on the study of the upper envelope of their future
infimum due to Pardo \cite{Pa}. These results extend integral test
and laws of the iterated logarithm for Bessel processes due to
Dvoretsky and Erd\"os \cite{de} and stable L\'evy processes
conditioned to stay positive with no positive
jumps due to Bertoin \cite{be1}.}\\


\noindent
{\it A.M.S. Classification}: {\footnotesize 60 G 18, 60 G 17, 60 G 51, 60 F 15.}\\

\author{J.C. Pardo.}
\section{Introduction.}A real Markov process $X=(X_{t}, t\geq 0)$ with
c\`adl\`ag paths is a self-similar process if for every $k>0$ and
every initial state $x\geq 0$ it satisfies the scaling property,
i.e., for some $\alpha>0$
\begin{equation*}
\textrm{the law of } (kX_{k^{-\alpha}t},t\geq 0) \textrm{ under
}\p_{x}\textrm{ is } \p_{kx},
\end{equation*}
where $\p_{x}$ denotes the law of the process $X$ starting from
$x\geq 0$.\\
In this article, we focus on positive self-similar Markov processes
and we will refer to them as {\it pssMp}. We denote by $X^{(x)}$ or
$(X, \p_x)$ for the pssMp starting from $x\geq 0$. Well-known
examples of pssMp  are: Bessel processes and stable subordinators or
more generally, stable processes conditioned to stay positive.\\

Self-similar processes were the object of a systematic study first
done by Lamperti~\cite{la1}. In a later work, Lamperti \cite{la}
studied the markovian case in detail. According to Lamperti
\cite{la} any pssMp starting from a strictly positive state
satisfies one of the following conditions:
\begin{itemize}
\item[i)] it never reaches the state $0$, \item[ii)] it hits the
state $0$ continuously or \item[iii)] it hits the  state $0$ by a
negative jump.
\end{itemize}
The main result in \cite{la} prove that any  pssMp starting from a
strictly positive state is a time-change of the exponential of a
L\'evy process. More precisely, let $X^{(x)}$ be a self-similar
Markov process started from $x>0$ that fulfills the scaling property
for some $\alpha>0$, then
\begin{equation}\label{lamp}
X^{(x)}_{t}=x\exp\Big\{\xi_{\tau(tx^{-\alpha})}\Big\},\qquad 0\leq
t\leq x^{\alpha}I(\xi),
\end{equation}
where
\begin{equation*}
\tau_{t}=\inf\Big\{s\geq 0: I_{s}(\xi)>t\Big\},\qquad
I_{s}(\xi)=\int_{0}^{s}\exp\big\{\alpha\xi_{u}\big\}\ud u,\qquad
I(\xi)=\lim_{t\to+\infty}I_{t}(\xi),
\end{equation*}
and $\xi$ is either a real L\'evy process which drifts towards
$-\infty$, if $X^{(x)}$ satisfies condition (ii) or $\xi$ is a
L\'evy process killed at an independent exponential time if
$X^{(x)}$ satisfies condition (iii) or $\xi$ is a L\'evy process
which does not drift towards $-\infty$, if $X^{(x)}$ satisfies
condition (i). This is the well-known
Lamperti representation.\\
Several authors have studied the problem of when an entrance law at
$0$ for $(X, \p_x)$ can be defined, see for instance Bertoin and
Caballero \cite{beC}, Bertoin and Yor \cite{BeY}, and Caballero and
Chaumont \cite{CCh}. Bertoin and Caballero \cite{beC} studied for
the first time this problem  for the increasing case. Later Bertoin
and Yor \cite{BeY} generalized the results obtained in \cite{beC}.
The main result of \cite{BeY} proves that this limiting process is
also a pssMp, and that it has the same semigroup as $(X, \p_x)$ for
$x>0$. Bertoin and Yor \cite{BeY} also gave sufficient conditions
for the weak convergence of $\p_{x}$ to hold when $x$ tends to $0$,
in the sense of finite dimensional distributions. The entrance law
was
also computed in the mentioned works, such law will be written below
 in (\ref{entrlaw}).\\
Caballero and Chaumont \cite{CCh} gave necessary and sufficient
conditions for the weak convergence of $X^{(x)}$ on the Skorokhod's
space. Caballero and Chaumont also gave a path construction of this
weak limit, that we will denote by $X^{(0)}$.\\

The aim of this work is to describe the upper envelope at $0$ and at
$+\infty$  for a large class of pssMp satisfying condition (i)
through integral test and laws of the iterated logarithm. We will
give special attention to the case with no positive jumps since  we
may obtain general integral tests and compare the rate of growth of
$X^{(x)}$ with that of its future infimum process and
the pssMp $X^{(x)}$ reflected at its future infimum.\\
Several partial results on the upper envelope of $X^{(0)}$ have
already been established before, the most important of which being
due to Dvoretsky and Erd\"os \cite{de} who studied the case of
Bessel process. More precisely, we have the well known Kolmogorov
and the Dvoretsky-Erd\"os integral test, see for instance It\^o and
McKean \cite{IMc}.
\begin{theorem}[Kolmogorov-Dvoretsky-Erd\"os] Let $h$ be a
nondecreasing, positive and unbounded function as $t$ goes to
$+\infty$. Then the upper envelope of $X^{(0)}$, a Bessel process
with index $\delta\ge 2$, at $0$ is as follows,
\[
\p\big(X^{(0)}_{t}>\sqrt{t}h(t), \textrm{ i.o., as } t\to 0\big)=0
\textrm{ or }1,
\]
according as,
\[
\int_{0}h^{\delta}(t)\exp\{-h^{2}(t)/2\}\frac{\ud t}{t}
\qquad\textrm{is finite or infinite.}
\]
\end{theorem}
Thanks to the time inversion property of Bessel processes, we also
have the integral test for large times. To get it, it is enough to
replace
\[
\int_{0}h^{\delta}(t)\exp\{-h^{2}(t)/2\}\frac{\ud t}{t}\quad
\textrm{ by }\quad
\int^{+\infty}h^{\delta}(t)\exp\{-h^{2}(t)/2\}\frac{\ud t}{t}.
\]
It is important to note that this integral test is also valid for
$d\geq 1$. The integral test for transient Bessel process will be
extended in section 4 (Theorems 2 and 3). Additionally in section 7,
we establish a variant of the
Kolmogorov-Dvoretsky-Erd\"os integral test, (Theorems 11 and 12).\\
From the  Kolmogorov-Dvoretsky-Erd\"os integral test, we deduce the
following laws of the iterated logarithm,
\begin{equation}\label{lawitlogbe}
\limsup_{t\to 0}\frac{X^{(0)}_{t}}{\sqrt{2t\log|\log
t|}}=1\quad\textrm{ and }\quad \limsup_{t\to
+\infty}\frac{X^{(0)}_{t}}{\sqrt{2t\log\log t}}=1,\quad \textrm{a.
s.}
\end{equation}

Recall that the future infimum of a pssMp starting at $x\geq 0$ is
defined by
\[
J^{(x)}_s =\inf_{t\ge s} X^{(x)}_t.
\]
Khohsnevisan, Lewis and Wembo \cite{Kh} studied the asymptotic
behaviour of the future imfimum of transient Bessel processes.
Khoshnevisan et al. \cite{Kh} established the following law of the
iterated logarithm for $J^{(0)}$ and for the Bessel process
reflected at its future infimum,  $X^{(0)}-J^{(0)}$;
\begin{equation}\label{lawitlogfibe}
\limsup_{t\to +\infty}\frac{J^{(0)}_{t}}{\sqrt{2t\log|\log
t|}}=1,\qquad \textrm{almost surely,}
\end{equation}
and
\begin{equation}\label{lawitlogdfibe}
\limsup_{t\to
+\infty}\frac{X^{(0)}_{t}-J^{(0)}_{t}}{\sqrt{2t\log|\log
t|}}=1,\qquad \textrm{almost surely.}
\end{equation}
In section 6, we extend these results and we also study the small
time behaviour.\\

When $X^{(0)}$ is a stable L\'evy process conditioned to stay
positive with no positive jumps and with index $1<\alpha\leq 2$, we
have the following law of the iterated logarithm due to Bertoin
\cite{be1},
\begin{equation}\label{lawitlogst}
\limsup_{t\to 0}\frac{X^{(0)}_{t}}{t^{1/\alpha}\big(\log|\log
t|\big)^{1-1/\alpha}}=c(\alpha)\qquad\textrm{ almost surely},
\end{equation}
where $c(\alpha)$ is a positive constant.\\
Recently in \cite{Pa} we studied the asymptotic behaviour of the
upper envelope of the future infimum of pssMp under general
hyphotesis. In particular, the author proved that when $X^{(0)}$ is
a stable L\'evy process conditioned to stay positive with no
positive jumps and with index $1<\alpha\leq 2$, we have
\begin{equation}\label{lawitlogfist}
\limsup_{t\to 0}\frac{J^{(0)}_{t}}{t^{1/\alpha}\big(\log|\log
t|\big)^{1-1/\alpha}}=c(\alpha),\qquad \textrm{almost surely,}
\end{equation}
and that
\begin{equation}\label{lawitlogfist2}
\limsup_{t\to +\infty}\frac{J^{(0)}_{t}}{t^{1/\alpha}\big(\log\log
t\big)^{1-1/\alpha}}=c(\alpha),\qquad \textrm{almost surely.}
\end{equation}
In section 6, we will see that under the assumption that
 $\p(J^{(0)}_{1}>t)$ is log-regular,i.e.
\[
-\log \p\big(J_1^{(0)}>t\big)\sim \lambda t^{\beta}L(t)\quad
\textrm{ as }\quad t\to +\infty,
\]
where $\lambda, \beta>0$ and $L$ is a function which varies slowly
at $+\infty$, the upper envelope of $X^{(x)}$ will be described by
an explicit law of the iterated logarithm and that it agrees with
the upper envelope of its future infimum as we have seen above.\\
All the asymptotic results presented in section 5 and 6 are
consequences of general integral tests which are stated in sections 3 and 4.\\
The rest of this paper is organized in six sections. Section 2 is
devoted to some preliminaries of L\'evy processes and pssMp. In
section 3, we study the asymptotic properties of the first and last
passage time processes of $X^{(0)}$. In section 4, we give the
general integral tests for the upper envelope of $X^{(x)}$. Sections
5 and 6 are devoted to applications of the results of sections 3 and
4. In section 7, we will apply our results to Bessel processes and
finally in section 8, we will give some examples.

\section{Preliminaries.}
\subsection{Weak convergence and entrance law of pssMp.}
Let $\mathcal{D}$ denote the space of Skorokhod of c\`adl\`ag paths.
We consider a probability measure on $\mathcal{D}$ denoted by $\p$
under which $\xi$ will always be a real L\'evy process such that
$\xi_0=0$. Let $\Pi$ be the L\'evy measure of $\xi$, that is the
measure satisfying
\[
\int_{(-\infty, \infty)}(1\land x^2)\Pi(\ud x)<\infty,
\]
and such that the characteristic exponent $\Psi$, defined by
\[
\e\Big(e^{iu\xi_t}\Big)= e^{-t\Psi(u)},\qquad t\geq 0, u\in \R
\]
is given, for some $b\geq 0$ and $a\in \R$, by
\[
\Psi(u)=iau+\frac 1 2 b^{2}u^2+\int_{(-\infty,
\infty)}\Big(1-e^{iux}+iux\ind_{\{|x|\leq 1\}}\Big)\Pi(\ud x), \quad
u\in \R.
\]
Define for $x\geq 0$,
\[
\overline{\Pi}^{+}(x)=\Pi\big((x,\infty)\big),\quad
\overline{\Pi}^{-}(x)=\Pi\big((-\infty,-x)\big), \quad
M(x)=\int_{0}^x \ud y \int_y^{\infty} \overline{\Pi}^-(z)\ud z,
\]
and
\[
J=\int_{[1,\infty)}\frac{x\overline{\Pi}^{+}(x)}{1+M(x)} \ud x.
\]
Then according to Caballero and Chaumont \cite{CCh}, necessary and
sufficient conditions for the weak convergence of $X^{(x)}$ on the
Skorokhod's space are
\[
\textrm{(\bf H)}\quad \xi\textrm{ is not arithmetic and }\quad
\begin{cases}
& \textrm{ either }0<\e(\xi_1)\leq \e(|\xi_1|)<\infty,\\
& \textrm{ or } \e(|\xi_1|)<\infty, \e(\xi_1)=0 \textrm{ and }
J<\infty,
\end{cases}
\]
and
\begin{equation}\label{condex}
\e\left(\log^+ \int_0^{T_1}\exp\big\{\alpha \xi_s\big\}\ud
s\right)<\infty,
\end{equation}
where $T_x$ is the first passage time above $x\geq 0$, i.e.
$T_x=\inf\{t\geq 0: \xi_t\geq x\}$.\\
The weak limit found in \cite{CCh}, denoted by $X^{(0)}$, is a pssMp
starting from $0$ which fulfills the Feller property on $[0,\infty)$
and with the same transition function as $X^{(x)}$, $x>0$. In all
the sequel, we suppose that the L\'evy processes considered here
satisfy conditions ({\bf H}) and (\ref{condex}). We will distinguish
the case $0<\e(\xi_1)\leq \e(|\xi_1|)<\infty$ saying that
$m:=\e(\xi_1)>0$.\\
It is important to note that if $m>0$, we have that
$\e(T_1)<\infty$ and hence condition (\ref{condex}) is satisfied.
Another example when this condition is satisfied is when $\xi$ has
no positive jumps (see section 2 in \cite{CCh}).\\

 We denote by $\p_{x}$ the law, under $\p$, of a
self-similar Markov process $X^{(x)}$ starting from $x>0$ and by
$\p_{0}$ the law, under $\p$, of the limiting process $X^{(0)}$.
When $m>0$ (i.e. $\xi$ has positive drift) the entrance law under
$\p_0$ has been computed by Bertoin and Caballero \cite{beC} and
Bertoin and Yor \cite{BeY} and can be expressed as follows: for
every $t>0$ and every measurable function $f:\R_{+}\to\R_{+}$,
\begin{equation}\label{entrlaw}
\e_{0}\big(f(X_{t})\big)=\frac1m
\e\Big(I^{-1}\big(\hat{\xi}\big)f\Big(tI^{-1}\big(\hat{\xi}\big)\Big)\Big)\qquad\textrm{for
}t>0,
\end{equation}
where
\[
I\big(\hat{\xi}\big)=\int_0^{+\infty}\exp\Big\{\alpha\hat{\xi}_s\Big\}\ud
s.
\]
When $m=0$ there is no explicit form for the entrance law of
$X^{(0)}$ in terms of the underlying L\'evy process. However,
Caballero and Chaumont \cite{CCh} proved that it can be obtained  as
the weak limit of the entrance law for the positive drift case, when
the drift tends towards 0. More precisely, for any bounded and
continuous function $f$,
\begin{equation}\label{entrlaw}
\e_{0}\big(f(X_{t})\big)=\lim_{\lambda \to 0}\frac1{\lambda}
\e\Big(I^{-1}_{\lambda}f\Big(tI^{-1}_{\lambda}\Big)\qquad\textrm{for
}t>0,
\end{equation}
where
\[
I_{\lambda}=\int_0^{+\infty}\exp\Big\{-\alpha(\xi_s-\lambda
s)\Big\}\ud s.
\]

Now, we introduce the so-called first and last passage times of
$X^{(0)}$ by
\[
S_y=\inf\Big\{t\geq 0: X^{(0)}_t\geq y\Big\}\quad\textrm{and}\quad
U_y=\sup\Big\{t\geq 0: X^{(0)}_t\leq y\Big\},
\]
for $y>0$. Note that when $m>0$, i.e. that the process $X^{(0)}$
drifts to $+\infty$, the last passage time $U_x$ is finite a.s. for
$x\geq 0$. Moreover, if the process $X^{(0)}$ satisfies the scaling
property with some index $\alpha>0$, then by the definition of $S_x$
and $U_x$, we deduce that the first passage time process $S=(S_x,
x\geq 0)$ and the last passage time process $U=(U_x, x\geq 0)$ are
increasing self-similar processes with scaling index $\alpha^{-1}$.
From the path properties of $X^{(0)}$ we easily see that both
processes start from $0$ and go to $+\infty$ as $x$
increases.\\
When $m=0$, the last passage time $U_x$ is no longer a.s. finite but
the first passage time process is still self-similar, starts from 0
and goes to $+\infty$ as $x$ increases.\\

{\it With no loss of generality, we will suppose that $\alpha=1$.
Indeed, we see from the scaling property that if $X^{(x)}$, $x\ge0$,
is a pssMp with index $\alpha>0$, then $\big(X^{(x)}\big)^{\alpha}$
is a pssMp with index equal to 1. Therefore, the integral tests and
LIL established in the sequel can easily be interpreted for any
$\alpha>0$. }

\subsection{PssMp with no positive jumps and some path transformations.}

Here, we suppose that $\xi$ has no positive jumps. It is known that
under this assumption, the process $\xi$ has finite exponential
moments of arbitrary positive order (see \cite{be} for background).
In particular, we have that
\begin{equation*}
\e\big(\exp\{u\xi_{t}\}\big)=\exp\{t\psi(u)\},\qquad u\geq 0
\end{equation*}
where $\psi$ is defined by,
\begin{equation*}
\psi(u)=au+\frac{1}{2}\sigma^{2}u^2
+\int_{]-\infty,0[}\big(e^{ux}-1-ux\ind_{\{x>-1\}}\big)\Pi(\ud x),
\quad u\geq 0.
\end{equation*}
It is important to note that the condition that $\xi$  does not
derive towards $-\infty$, is equiva\-lent to
\begin{equation*}
m=\psi'(0+)\in [0,\infty[.
\end{equation*}
We recall  Caballero and Chaumont's construction  of $X^{(0)}$, only
in this particular case.  In this direction, let $T_{z}=\inf\{s\geq0
: \xi_{s}\geq z\}$ be the first time where the process $\xi$ reaches
the state $z\geq 0$. Note that due to the absence of positive jumps
and since the process $\xi$ does not derive towards $-\infty$, then
for all $z\geq 0$
\begin{equation*}
T_{z}<\infty \quad \textrm{and}\quad\xi_{T_{z}}=z,\quad
\p-\textrm{a.s.}
\end{equation*}
Let $x_{1}\geq x_{2}\geq \cdots>0$, be an infinite decreasing
sequence of strictly positive real numbers which converges to $0$
and $(\xi^{(n)}, n\geq 1)$ a sequence of random processes  which
are independent and have the same distribution as $\xi$. \\
From the sequences $(x_{n})$ and $\big(\xi^{(n)}, n\geq 1\big)$, and
by Lamperti's transformation (\ref{lamp}) we may define a sequence
of pssMp. More precisely,
\begin{equation}\label{lampp}
X^{(x_{n})}_{t}=x_{n}\exp\Big\{\xi^{(n)}_{\tau^{(n)}(t/x_{n})}\Big\},
\quad t\geq 0, \quad, n\geq 1,
\end{equation}
where
\begin{equation*}
\tau^{(n)}(t)=\inf\Big\{s\geq 0: I_{s}(\xi^{(n)})>t\Big\}, \qquad
I_{t}(\xi^{(n)})=\int_{0}^t\exp\Big\{ \xi^{(n)}_{u}\Big\}\ud u.
\end{equation*}
For each $n\geq 2$, we also define the first time in which the
process $X^{(x_{n})}$ reaches the state $x_{n-1}$, i.e.
\begin{equation*}
S^{(n)}=\inf\Big\{t\geq 0:X^{(x_{n})}_{t}\geq x_{n-1}\Big\}.
\end{equation*}
It is clear from (\ref{lampp}) that
$S^{(n)}=x_{n}I_{T^{(n)}}(\xi^{(n)})$, where $T^{(n)}$ is the first
passage time where the process $\xi^{(n)}$ reaches the state
$\log(x_{n-1}/x_{n})$, i.e.
\[
T^{(n)}=\inf\{t\geq 0:\xi^{(n)}\ge \log(x_{n-1}/x_{n})\}.
\]
Now from the sequences  $\big(X^{(x_{n})}, n\geq 1\big)$ and
$\big(S^{(n)}, n\geq 2\big)$, we can construct the process $X^{(0)}$
as the concatenation of the processes $X^{(x_{n})}$ on each interval
$\big[0,S^{(n)}\big]$, i.e.,
\begin{equation*}
X_{t}^{(0)}=\left\{\begin{array}{ll}
X_{t-\Sigma_{2}}^{(x_{1})} & \textrm{ if } t\in[\Sigma_{2},\infty[,\\
X_{t-\Sigma_{3}}^{(x_{2})} & \textrm{ if }
t\in[\Sigma_{3},\Sigma_{2}[,\\\vdots & \\
X_{t-\Sigma_{n+1}}^{(x_{n})} & \textrm{ if }
t\in[\Sigma_{n+1},\Sigma_{n}[,\\
\vdots &
\end{array} \right.
\end{equation*}
where $\Sigma_{n}=\sum_{k\geq n}S^{(k)}$. We can deduce, by the
definition of the process $X^{(0)}$, that
$\Sigma_{n}=\inf\big\{t\geq 0: X^{(0)}_{t}\geq
x_{n-1}\big\}$\\
Caballero and Chaumont proved that  this construction  makes sense,
it does not depend on the sequence $(x_{n})$, and $X^{(0)}_{0}=0$.
As we mentioned before,  $X^{(0)}$ is a c\`adl\`ag pssMp with the
same semi-group as $(X,\p_{x})$ for $x>0$. In the general case such
construction is more complicate since it depends also on the
overshoots of the sequence of L\'evy processes $(\xi^{(n)}, n\geq 1)$.\\
Note that the process $X^{(0)}$ inherits the path properties of the
underlying L\'evy process, hence $X^{(0)}$ does not have positive
jumps and if $m>0$ it drifts towards $+\infty$. Hence for all $x\geq
0$
\[
S_x\textrm{ and } U_x\textrm{ are finite}\quad\textrm{and}\quad
X^{(0)}_{S_x}=X^{(0)}_{U_x}=x, \quad\textrm{almost surely.}
\]
From this construction and the path decomposition in Corollary 1 in
Chaumont and Pardo \cite{CP}, we deduce that the first passage time
process $S=(S_x, x\geq 0)$ and the last passage time process
$U=(U_x, x\geq 0)$ are increasing self-similar processes with
independent increments, in fact we will prove at the end of this
section that these processes are also self-decomposable. This
property was studied by the first time by Getoor in \cite{Get} for
the last passage time of a Bessel process of index $\delta\geq 3$
and later by Jeanblanc, Pitman and Yor \cite{JPY} for $\delta>2$. If
$m=0$, the process $S$ is still an increasing self-similar process
with
independent increments and also self-decomposable.\\

We are interested in describing the law of the process
($X^{(0)}_{(S_x-t)^-}, 0\leq t\leq S_x$) and in obtaining the law of
the first passage time in terms of the underlying L\'evy process.
With this purpose, we now briefly recall the definition of the
L\'evy process conditioned to stay positive $\xi^\uparrow$ and refer
to \cite{CD} for a complete account on this subject. The process
$\xi^\uparrow$ is an $h$-process of $\xi$ killed when it first
exists $(0,\infty)$, i.e. at time $R=\inf\{t:\xi_t\le0\}$. The law
of this strong Markov process $\xi^\uparrow$ is defined by its
semi-group:
\[\p(\xi^\uparrow_{t+s}\in dy\,|\xi^\uparrow_s=x)=
\frac{h(y)}{h(x)}\p(\xi_t+x\in dy\,,t<R)
\,,\;\;\;s,t\ge0,\,\;\;\;x,y>0\] and its entrance law:
\[\p(\xi^\uparrow_t\in dx)=h(x)\hat{N}(\xi_t\in dx,\,t< \zeta)\,,\]
where $\hat{N}$ is the excursion measure of the reflected process
$\xi-I=(\xi_t-\inf_{s\leq t}\xi_s, t\geq 0)$, $\zeta$ is the
lifetime of the generic excursion and $h$ is the positive harmonic
function (for $\xi$ killed at time $R$) which is defined by:
\[h(x)=\e\left(\int_0^\infty\ind_{\{I_t\ge -x\}}\,dL_t\right)
\,,\;\;\;x\ge0\, ,
\]
where $I_t= \inf_{s\leq t}\xi_s$ and $L$ is the local time of the
reflected process $\xi-I$.\\
Note that $\xi^\uparrow$ has no positive jumps and that almost
surely
\[
\lim_{t\downarrow0}\xi^\uparrow_t=0,\qquad
\lim_{t\uparrow+\infty}\xi_t^\uparrow=+\infty,\quad \textrm{ and }
\quad \xi^\uparrow_t>0\quad \textrm{ for all }t>0.
\]
The following time reversal property of $\xi^\uparrow$ is an
important tool for our next result, its proof can be found in
Theorem VII.18 in Bertoin \cite{be}
\begin{lemma}
The law of $(x-\xi_t, 0\leq t\leq T_x )$ is the same as that of the
time-reversed process $(\xi^\uparrow_{(\gamma^\uparrow(x)-t)-},
0\leq t \leq \gamma^\uparrow(x))$, where
$\gamma^\uparrow(x)=\sup\{t\geq 0, \xi^\uparrow_{t}\leq x\}$.
Moreover, for every $x>0$, the process $\big(\xi^\uparrow_{t}, 0\leq
t\leq \gamma^\uparrow(x)\big)$ is independent of
$\big(\xi^\uparrow_{\gamma^\uparrow(x)+t}-x, t\geq 0\big)$, and the
latter has the same law as that of the process $\xi^\uparrow$.
\end{lemma}
Now, for every $y>0$ let us define
\[
\tilde{X}^{(y)}_t=y\exp\Big\{-\xi^\uparrow_{\tau^\uparrow(t/y)}\Big\}\qquad
t\geq 0,
\]
where
\begin{equation*}
\tau^\uparrow_{t}=\inf\Big\{s\geq 0:
I_{s}\big(-\xi^\uparrow\big)>t\Big\},\quad \textrm{ and }\quad
I_{s}\big(-\xi^\uparrow\big)=\int_{0}^{s}\exp\big\{-\xi^\uparrow_{u}\big\}\ud
u.
\end{equation*}
We denote by $\tilde{\p}_y$, the law of $\tilde{X}^{(y)}$. Since
$\xi^\uparrow$ derives towards $+\infty$, we deduce that
$\tilde{X}^{(y)}$, reaches $0$ at an almost surely finite random
time, denoted by $\tilde{\rho}_y=\inf\{t\geq 0,
\tilde{X}^{(y)}_{t}=0\}$.
\begin{proposition}
Suppose $m\geq 0$. The law of the process time-reversed at its first
passage time below x, $\big(X^{(0)}_{(S_x-t)-}, 0\leq t\leq
S_x\big)$ is the same as that of the process $(\tilde{X}^{(x)}_{t},
0\leq t\leq \tilde{\rho}_x)$.
\end{proposition}
\noindent\textit{Proof: } Let us take any decreasing sequence
$(x_{n})$ of positive real numbers which converges to $0$ and such
that $x_{1}=x$. By Lemma 1, we can split $(\tilde{X}^{(x)}_{t},
0\leq t\leq \tilde{\rho})$ into the sequence
\begin{equation*}
\Big(x_{1}\exp\big\{-\xi^\uparrow_{\tau^\uparrow(t/x_{1})}\big\},
x_{1}I_{\gamma^{(n)}}\big(-\xi^\uparrow\big)\leq t \leq
x_{1}I_{\gamma^{(n+1)}}\big(-\xi^\uparrow\big)\Big), \quad n\geq 1,
\end{equation*}
where $\gamma^{(n)}=\sup\{t\geq0: \xi^\uparrow_{t}\leq \log
x_{n}/x\}$.\\
Then to prove this result, it is enough to show that, for each
$n\geq 1$
\begin{equation*}
\Big(X^{(0)}_{(S_{x_{n}}-t)-}, 0\leq t\leq S^{(n)}\Big)
\overset{(d)}{=}\Big(x\exp\big\{-\xi^\uparrow_{\tau^\uparrow(t/x_{1})}\big\},
xI_{\gamma^{(n)}}\big(-\xi^\uparrow\big)\leq t \leq
xI_{\gamma^{(n+1)}}\big(-\xi^\uparrow\big)\Big)
\end{equation*}
where $S^{(n)}=S_{x_{n}}-S_{x_{n+1}}$.\\
Fix $n\geq 1$, from the Caballero and Chaumont's construction, we
know that the left-hand side of the above identity has the same law
as
\begin{equation}\label{const3}
\bigg(x_{n+1}\exp\Big\{{\xi}^{(n+1)}_{\tau^{(n+1)}\big(I_{T^{(n+1)}}\big(\xi^{(n+1)}\big)-t/x_{n+1}\big)}\Big\},
0 \leq t \leq x_{n+1}I_{T^{(n+1)}}\big(\xi^{(n+1)}\big)\bigg).
\end{equation}
On the other hand by Lemma 1, we know that $(-\xi^\uparrow_{t},
0\leq t \leq \gamma^{(n)})$ is independent of
$\xi^{\uparrow(n)}=(\log (x/x_{n})-\xi^\uparrow_{\gamma^{(n)}+t},
t\geq 0)$ and that the latter has the same law as $-\xi^\uparrow$.
Since
\[
\tau^\uparrow\Big(I_{\gamma^{(n)}}\big(-\xi^\uparrow\big)+t/x\Big)=\gamma^{(n)}+\inf\left\{s\geq
0: \int_{0}^{s}\exp\Big\{-\xi^{^\uparrow(n)}_{u}\Big\}\ud u\geq
t/x_{n}\right\},
\]
it is clear that the right-hand side of the above identity in
distribution has the same law as,
\begin{equation}\label{ddual1}
\Big(x_{n}\exp\big\{-\xi^\uparrow_{\tau^\uparrow(t/x_{n})}\big\},
0\leq t\leq
x_{n}I_{\gamma^\uparrow(\log(x_{n}/x_{n+1}))}\big(-\xi^\uparrow\big)\Big).
\end{equation}
Therefore, it is enough to show that (\ref{const3}) and
(\ref{ddual1}) have the same distribution.\\
Now, let us define the exponential functional of
$\big(\xi^{(n+1)}_{(T^{(n+1)}-t)^-}, 0\leq t\leq T^{(n+1)}\big)$ as
follows,
\begin{equation*}
B^{(n+1)}_{s}=\int_{0}^{s}\exp\Big\{\xi^{(n+1)}_{T^{(n+1)}-u}\Big\}\ud
u\qquad \textrm{for } s\in[0,T^{(n+1)}],
\end{equation*}
and $\mathcal{H}(t)=\inf\big\{ 0\leq s\leq T^{(n+1)},
B^{(n+1)}_{s}>t\big\}$, the right continuous inverse of the
exponential functional $B^{(n+1)}$.\\
By a change of variable, it is clear that
$B^{(n+1)}_{s}=I_{T^{(n+1)}}\big(\xi^{(n+1)}\big)-I_{T^{(n+1)}-s}\big(\xi^{(n+1)}\big),$
and if we set $t=x_{n+1}B^{(n+1)}_{s}$, then
$s=\mathcal{H}(t/x_{n+1})$ and hence
\begin{equation*}
\begin{split}
\tau^{(n+1)}\Big(I_{T^{(n+1)}}\big(\xi^{(n+1)}\big)-t/x_{n+1}\Big)&
=\tau^{(n+1)}\Big(I_{T^{(n+1)}-s}\big(\xi^{(n+1)}\big)\Big)\\
&=T^{(n+1)}-\mathcal{H}(t/x_{n+1}).
\end{split}
\end{equation*}
Therefore, we can rewrite (\ref{const3}) as follows
\begin{equation}\label{const3rt}
\bigg(x_{n+1}\exp\Big\{\xi^{(n+1)}_{T^{(n+1)}-\mathcal{H}(t/x_{n+1})}\Big\},
0 \leq t \leq x_{n+1}B^{(n+1)}_{T^{(n+1)}}\bigg),
\end{equation}
and applying Lemma 1, we get that (\ref{const3rt}) has the same law
as that of the process defined in (\ref{ddual1}).\QED \noindent An
important consequence of this Proposition is the following
time-reversed identity. For any $y < x$,
\begin{equation*}
\Big(X^{(0)}_{(S_{x}-t)-}, S_{y}\leq t\leq
S_{x}\Big)\overset{(d)}{=}\big(\tilde{X}^{(x)}_{t},0\leq t\leq
\tilde{U}_{y}\big),
\end{equation*}
where $\tilde{U}_{y}=\sup\big\{t\geq 0,\tilde{X}^{(x)}_{t}\leq y
\big\}$.\\
Four the next result we need to introduce the concept of
self-decomposable random variable.
\begin{defi}
We say that a random variable $X$ is self-decomposable if for every
constant $0<c<1$ there exists a variable $Y_{c}$ which is
independent of $X$ and such that $Y_{c}+cX$ has the same law as $X$.
\end{defi}

\begin{corollary}Let $m\geq 0$. For every $x>0$ the law  of $S_x$, the first passage time
of the process $X^{(0)}$ above $x$, has the same law as
$xI\big(-\xi^\uparrow\big)$, where
\begin{equation*}
I\big(-\xi^\uparrow\big)=\int_{0}^{\infty}\exp\{-\xi^\uparrow_{u}\}\ud
u.
\end{equation*}
\end{corollary}
\noindent\textit{Proof: }From Proposition 1, we see that $S_x$ and
$\tilde{\rho}_x$ have  the same law. By the Lamperti representation
of $\tilde{X}^{(x)}$, we deduce that
$\tilde{\rho}_x=xI\big(-\xi^\uparrow\big)$ and then the identity in
law
follows.\\
Now, let $0<c<1$. From Lemma 1, we know that the killed process
$(\xi^\uparrow_{t}, 0\leq t\leq \gamma^\uparrow(\log(1/c)))$ is
independent of $(\xi^\uparrow_{\gamma^\uparrow(\log(1/c))+t}+\log c,
t\geq 0)$ and that the latter has the same as $\xi^\uparrow$, then
\[
I\big(-\xi^\uparrow\big)=\int_{0}^{\gamma^\uparrow(\log(1/c))}\exp\{-\xi^\uparrow_{u}\}\ud
u+c\int_{0}^{+\infty}\exp\big\{-\xi^\uparrow_{\gamma^\uparrow(\log(1/c))+u}-\log
c\big\}\ud u,
\]
the self-decomposability follows.\QED \noindent
\begin{lemma}Let $m>0$. For every $x>0$ the law  of $U_x$, the last passage time of
the process $X^{(0)}$ below $x$, has the same law as
$xI\big(\hat{\xi}\big)$. Moreover, $U_1$ is self-decomposable.
\end{lemma}
\noindent\textit{Proof: }The first part of this Lemma is consequence
of Proposition 1 in \cite{CP}. Let $0<c<1$. From the Markov
property, we know that $(\xi_{t}, 0\leq t\leq T_{\log (1/c)})$ is
independent of $(\xi_{T_{\log (1/c)}+t}+\log c, t\geq 0)$ and that
the latter has the same distribution as $\xi$, then
\[
I\big(\hat{\xi}\big)=\int_{0}^{T_{\log(1/c)}}\exp\big\{-\xi_{u}\big\}\ud
u +c\int_{0}^{+\infty}\exp\big\{-\xi_{T_{\log (1/c)}+u}-\log
c\big\}\ud u,
\]
the self-decomposability follows.\QED
\section{The lower envelope of the first passage time.}
In this section, we are interested in describing the lower envelope
at $0$ and at $+\infty$ of the first passage time of $X^{(0)}$
trough integral tests. We first deal with the case when $X^{(0)}$
has no positive jumps since in this case we may obtain an explicit
integral tests in terms of the tail probability  of
$I(-\xi^{\uparrow})$ and also we may consider the case when $m=0$.
The general case will not be developed in a complete form but we
will use the integral tests of the lower envelope of the last
passage times due to Pardo \cite{Pa} (see Theorems 3 and 4) to get
the upper bound. In fact, it does not seem easy to determine the law
of $S_1$ in terms of the underlying L\'evy process and even
establish integral tests in terms of the decomposition of the first
passage time  in the Caballero and Chaumont's construction since the
sequence of the overshoots related to the underlying L\'evy process
is a Markov chain (see Proposition 2 in \cite{CCh}). This limited to
us to assume that $m>0$, since the last passage time is not defined
on the oscillating case.
\subsection{The case with no positive jumps.}
Let us define
\[
F^\uparrow(t)\eqdef\p\left(I(-\xi^\uparrow)<t\right),
\] and we denote by
$\mathcal{H}^{-1}_{0}$, the totality of positive increasing
functions $h(x)$ on $(0, \infty)$ that satisfy
\begin{itemize}
\item[i)] $h(0)=0$, and \item[ii)] there exists $\beta\in(0,1)$
such that
$\displaystyle\sup_{t<\beta}\displaystyle\frac{h(t)}{t}<\infty.$
\end{itemize}
The lower envelope at $0$ of the first passage process $S$ is as
follows:
\begin{proposition}
Let $m\geq 0$ and $h\in\mathcal{H}^{-1}_{0}$.
\begin{itemize}
\item[i)] If
\[
\int_{0^{+}}F^\uparrow\left(\frac{h(x)}{x}\right)\frac{\ud
x}{x}<\infty,
\]
then for all $\epsilon > 0$
\[
\p\Big(S_x<(1-\epsilon)h(x),\textrm{ i.o., as }x\to 0\Big)=0.
\]
\item[ii)] If
\[
\int_{0^{+}}F^\uparrow\left(\frac{h(x)}{x}\right)\frac{\ud
x}{x}=\infty,
\]
then for all $\epsilon > 0$
\[
\p\Big(S_x<(1+\epsilon)h(x),\textrm{ i.o., as }x\to 0\Big)=1.
\]
\end{itemize}
\end{proposition}
Let us define $\mathcal{H}^{-1}_{\infty}$, the totality of positive
increasing functions $h(t)$ on $(0, \infty)$ that satisfy
\begin{itemize}
\item[i)] $\lim_{t\to \infty}h(t)=+\infty$, and \item[ii)] there
exists $\beta\in(1,+\infty)$ such that
$\displaystyle\sup_{t>\beta}\displaystyle\frac{h(t)}{t}<\infty.$
\end{itemize}
The  lower envelope at $+\infty$ of the first passage process $S$ is
as follows:
\begin{proposition}
Let $m\geq 0$ and $h\in\mathcal{H}^{-1}_{\infty}$.
\begin{itemize}
\item[i)] If
\[
\int^{+\infty}F^\uparrow\left(\frac{h(x)}{x}\right)\frac{\ud
x}{x}<\infty,
\]
then for all $\epsilon > 0$
\[
\p\Big(S_x<(1-\epsilon)h(x),\textrm{ i.o., as }x\to +\infty\Big)=0.
\]
\item[ii)] If
\[
\int_{+\infty}F^\uparrow\left(\frac{h(x)}{x}\right)\frac{\ud
x}{x}=\infty,
\]
then for all $\epsilon > 0$
\[
\p\Big(S_x<(1+\epsilon)h(x),\textrm{ i.o., as }x\to +\infty\Big)=1.
\]
\end{itemize}
\end{proposition}

The above Propositions are consequence of Lemmas 3.1 and 3.2 of
Watanabe \cite{wa} and Corollary 1, this follows from the fact that
$S$ is an increasing self-similar process. It is important to note
that the above results may be proved using similar arguments as in
Theorems 3 and 4 in Pardo \cite{Pa}, it is enough to exchange
$I(\hat{\xi})$ by $I(-\xi^\uparrow)$ and note that $\Gamma=1$.
\subsection{The general case.}
Let us define
\[
G(t)\eqdef \p\big( S_{1} < t\big)\quad\textrm{ and } \quad
\overline{F}(t)\eqdef \p\big( I(\hat{\xi}) < t\big).
\]
The  lower envelope at $0$ of the first passage process $S$ is as
follows:
\begin{proposition}  Let $m > 0$ and
$h\in\mathcal{H}^{-1}_{0}$.
\begin{itemize}
\item[i)] If
\[
\int_{0^{+}}G\left(\frac{h(x)}{x}\right)\frac{\ud x}{x}<\infty,
\]
then for all $\epsilon > 0$
\[
\p\Big(S_{x}<(1-\epsilon)h(x),\textrm{ i.o., as } x\to 0\Big)=0.
\]
\item[ii)]If
\[
\int_{0^{+}}\overline{F}\left(\frac{h(x)}{x}\right)\frac{\ud
x}{x}=\infty,
\]
then for all $\epsilon > 0$
\[
\p\Big(S_x<(1+\epsilon)h(x),\textrm{ i.o., as }x\to 0\Big)=1.
\]
\end{itemize}
\end{proposition}
\textit{Proof:} We first prove the convergent part. Let $(x_{n})$ be
a decreasing sequence of positive numbers which converges to $0$ and
let us define the
events $A_{n}=\{S_{x_{n+1}}<h(x_{n})\}$.\\
Now, we choose $x_{n}=r^{n}$, for $r<1$. From the first Borel
Cantelli's Lemma, if $\sum_{n}\p(A_{n})<\infty$, it follows
\[
S_{r^{n+1}}\geq h\big(r^{n}\big)\qquad \p-\textrm{ a.s.,}
\]
for all large $n$. Since the function $h$ and the process $S$ are
increasing, we have
\[
S_x\geq h(x)\qquad \textrm{for}\quad r^{n+1}\leq x\leq r^{n}.
\]
On the other hand, from the scaling property, we get that
\[
\begin{split}
\sum_{n}\p\Big(S_{r^{n}}<h\big(r^{n+1}\big)\Big)&\leq
\int_{1}^{\infty}\p\Big(r^{t}S_{1}<h\big(r^{t}\big)\Big)\ud t\\
&=-\frac{1}{\log{r}}\int^{r}_{0}G\left(\frac{h(x)}{x}\right)\frac{\ud
x}{x}.
\end{split}
\]
From our hypothesis, this last integral is finite. Then from the
above discussion, there exist $x_{0}$ such that for every $x\geq
x_{0}$
\[
S_x\geq r^{2}h(x),\qquad\textrm{for all }\quad r<1.
\]
Clearly, this implies that
\[
\p_{0}\Big( S_x<r^{2}h(x), \textrm{ i.o., as }x\to 0\Big)=0,
\]
which proves part $(i)$.\\
The divergent part follows from the integral test for the lower
envelope of the last passage time due to Pardo \cite{Pa}, ( see
Theorem 3, part $(ii)$).\QED

The  lower envelope at $+\infty$ of the first passage process $S$ is
as follows:
\begin{proposition}  Let $m>0$ and
$h\in\mathcal{H}^{-1}_{\infty}$.
\begin{itemize}
\item[i)] If
\[
\int^{+\infty}G\left(\frac{h(x)}{x}\right)\frac{\ud x}{x}<\infty,
\]
then for all $\epsilon > 0$
\[
\p\Big(S_{x}<(1-\epsilon)h(x),\textrm{ i.o., as }x\to
+\infty\Big)=0.
\]
\item[ii)]If
\[
\int^{+\infty}\overline{F}\left(\frac{h(x)}{x}\right)\frac{\ud
x}{x}=\infty,
\]
then for all $\epsilon > 0$
\[
\p\Big(S_x<(1+\epsilon)h(x),\textrm{ i.o., as }x\to +\infty\Big)=1.
\]
\end{itemize}
\end{proposition}
\textit{Proof:} The proof is very similar to that in Proposition 4.
We get the integral test following the same arguments for the proof
of part $(i)$ and $(ii)$ for the sequence $x_{n}=r^{n}$, with $r>1$.
\QED  Note that in the general case, the integral tests for the
lower envelope of $S$ no longer depend on $F^\uparrow$ as in the
case with no positive jumps. Recall that it does not seem easy to
determine the law of $S$ and even have a nice decomposition as for
the last passage times that allows us to obtain
and integral test which only depends on the law of $S_1$.\\
We remark that $\overline{F}(t)$ is smaller or equal to $G(t)$, for
$t\geq 0$, but for our purpose the integral tests of above will be
very useful since we will compare in sections 5 and 6 the behaviour
of $\overline{F}$ and $G$ under different conditions.
\section{The upper envelope of pssMp.}
Here, we are interested in describing the upper envelope at $0$ and
at $+\infty$ of the pssMp $X^{(x)}$ trough integral tests. By the
same reasons as those mentioned in the precedent section, we will
first study the case with no positive jumps.
\subsection{The case with no positive jumps.}
The following theorem means in particular that the upper envelope at
$0$ of $X^{(0)}$ only depends on the tail behaviour of the law of
$I(-\xi^\uparrow)$ and on the additional hypothesis
\begin{equation}\label{adhyp}
\e\Big(\log^{+}I\big(-\xi^\uparrow\big)^{-1}\Big)<\infty.
\end{equation}
Let us recall that
\[
F^\uparrow(t)=\p\Big( I\big(-\xi^\uparrow\big)< t\Big),
\]
and denote by $\mathcal{H}_{0}$ the totality of positive increasing
functions $h(t)$ on $(0, \infty)$ that satisfy
\begin{itemize}
\item[i)] $h(0)=0$, and \item[ii)] there exists $\beta\in(0,1)$
such that
$\displaystyle\sup_{t<\beta}\displaystyle\frac{t}{h(t)}<\infty.$
\end{itemize}
\newpage
\begin{theorem}
Let $m\geq 0$ and $h\in\mathcal{H}_{0}$.
\begin{itemize}
\item[i)] If
\[
\int_{0^{+}}F^\uparrow\left(\frac{t}{h(t)}\right)\frac{\ud
t}{t}<\infty,
\]
then for all $\epsilon > 0$
\[
\p_{0}\Big(X_{t}>(1+\epsilon)h(t),\textrm{ i.o., as }t\to 0\Big)=0.
\]
\item[ii)] Assume that (\ref{adhyp}) is satisfied. If
\[
\int_{0^{+}}F^\uparrow\left(\frac{t}{h(t)}\right)\frac{\ud
t}{t}=\infty,
\]
then for all $\epsilon > 0$
\[
\p_{0}\Big(X_{t}>(1-\epsilon)h(t),\textrm{ i.o., as }t\to 0\Big)=1.
\]
\end{itemize}
\end{theorem}
\textit{Proof:} Let $(x_{n})$ be a decreasing sequence which
converges to $0$. We define the events $A_{n}=\big\{\textrm{There
exists }t\in [S_{x_{n+1}}, S_{x_{n}}]\textrm{ such that }
X^{(0)}_{t}>h(t)\big\}$. From the fact that $S_{x_{n}}$ tends to
$0$, a.s. when $n$ goes to $+\infty$, we see
\[
\Big\{X^{(0)}_{t}>h(t), \textrm{ i.o., as }t\to
0\Big\}=\limsup_{n\to +\infty} A_{n}.
\]
Since $h$ is an increasing function the following inclusions hold
\begin{equation}\label{inclus}
\Big\{x_{n}>h\big(S_{x_{n}}\big) \Big\}\subset A_{n}\subset
\Big\{x_{n}>h\big(S_{x_{n+1}}\big)\Big\}.
\end{equation}
Now, we prove the convergent part. We choose $x_{n}=r^{n}$, for
$r<1$ and $h_{r}(t)=r^{-2}h(t)$. Since $h$ is increasing, we deduce
that
\[
\sum_{n}
\p\Big(r^n>h_{r}\big(S_{r^{n+1}}\big)\Big)\leq\int_{1}^{+\infty}\p\Big(r^t>h\big(S_{t^r}\big)\Big)\ud
t\leq -\frac{1}{\log
r}\int_{0}^{r}\p\Big(t>h\big(S_t\big)\Big)\frac{\ud t}{t}.
\]
Replacing $h$ by $h_r$ in (\ref{inclus}), we see that we can obtain
our result if
\[
 \int_{0}^{r}\p\Big(t>h\big(S_t\big)\Big)\frac{\ud t}{t}<\infty.
\]
From elementary calculations and Corollary 1, we deduce that
\[
\int_{0}^{r}\p\Big(t>h\big(S_t\big)\Big)\frac{\ud t}{t}=\e \left(
\int_{0}^{h^{-1}(r)}\ind_{\big\{t/r<
I(-\xi^\uparrow)<t/h(t)\big\}}\frac{\ud t}{t}\right),
\]
where $h^{-1}(s)=\inf\{t>0, h(t)>s\}$, the right inverse function of
$h$. Then, this integral converges if
\[
\int_{0}^{h^{-1}(r)}\p\left(
I\big(-\xi^\uparrow\big)<\frac{t}{h(t)}\right)\frac{\ud t}{t}
<\infty.
\]
This proves part $(i)$.\\
Next, we prove the divergent case. We suppose that $h$ satisfies
\[
\int_{0^+}F^\uparrow\left(\frac{t}{h(t)}\right)\frac{\ud
t}{t}=\infty.
\]
Take, again, $x_{n}=r^{n}$, for $r<1$ and  define
\[
B_{n}\eqdef\bigcup_{m=n}^{\infty}A_{m}=\Big\{\textrm{There exists
}t\in (0, S_{r^n}]\textrm{ such that } X^{(0)}_{t}>h_{r}(t)\Big\}.
\]
Note that the family $(B_n)$ is decreasing and
\[
B\eqdef\bigcap_{n\geq 1}B_{n}=\Big\{X^{(0)}_{t}>h_{r}(t), \textrm{
i.o., as }t\to 0\Big\},
\]
then it is enough to prove that $\lim\p(B_{n})=1$ to obtain our
result.\\
Again replacing $h$ by $h_{r}$ in inclusion (\ref{inclus}), we see
\begin{equation}\label{desteo}
\p(B_{n})\geq 1-\p\Big(r^{j}\leq h_{r}\big(S_{r^j})\big), \textrm{
for all } n\leq j\leq m-1 \Big),
\end{equation}
where $m$ is chosen arbitrarily $m\geq n+1$.\\
Now, we define the events
\[
C_{n}\eqdef\bigg\{r^n>rh\Big(S_{r^{n}})\Big)\bigg\}.
\]
We will prove that $\sum \p(C_{n})=\infty.$  Since the function $h$
is increasing, from the identity in law of Corollary 1 it is
straightforward that
\[
\sum_{n\geq 1} \p(C_{n}) \geq
\int_{0}^{+\infty}\p\left(r^{t}>h\Big(S_{r^{t}}\big)\Big)\right)\ud
t=-\frac{1}{\log
r}\int_{0}^{1}\p\Big(t>h\big(tI\big(-\xi^\uparrow\big)\big)\Big)\frac{\ud
t}{t}.
\]
Hence, if this last integral is infinite, we get that $\sum
\p(C_{n})=\infty.$ In this direction, we have
\[
\int_{0}^{r}\p\Big(t>h\big(tI\big(-\xi^\uparrow\big)\big)\Big)\frac{\ud
t}{t} =\e \left( \int_{0}^{h^{-1}(r)}\ind_{\big\{t/r<
I(-\xi^\uparrow)<t/h(t)\big\}}\frac{\ud t}{t}\right).
\]
On the other hand, we see
\[
\begin{split}
\int_{0}^{h^{-1}(r)}\p\left(
I\big(-\xi^\uparrow\big)<\frac{t}{h(t)}\right)\frac{\ud t}{t}=
&\int_{0}^{h^{-1}(r)}\p\left(\frac{t}{r}<
I\big(-\xi^\uparrow\big)<\frac{t}{h(t)}\right)\frac{\ud
t}{t}\\&+\int_{0}^{h^{-1}(r)}\p\left(
I\big(-\xi^\uparrow\big)<\frac{t}{r}\right)\frac{\ud t}{t},
\end{split}
\]
and
\[
\int_{0}^{h^{-1}(r)}\p\left(
I\big(-\xi^\uparrow\big)<\frac{t}{r}\right)\frac{\ud t}{t}\leq
\e\left(\log^{+}
\frac{h^{-1}(r)}{r}I\big(-\xi^\uparrow\big)^{-1}\right)
\]
which is clearly finite from our assumptions. Then, we deduce
\[
\e \left( \int_{0}^{h^{-1}(r)}\ind_{\big\{t/r<
I(-\xi^\uparrow)<t/h(t)\big\}}\frac{\ud t}{t}\right)=\infty,
\]
and hence $\sum\p(C_{n})=\infty$.\\
Next, for $n\leq m-1$, we define
\[
H(n,m)\eqdef\p\Big(r^{j}\leq rh\big(S_{r^{j}}-S_{r^{m}}\big),
\textrm{ for all }n\leq j\leq m-1\Big),
\]
and we will prove that there exist $(n_l)$ and $(m_{l})$, two
increasing sequences such that $0\leq n_{l}\leq m_{l}-1$, and
$n_{l}, m_{l}$ go to $+\infty$ and
$H(n_l,m_l)$ tends to $0$ as $l$ goes to infinity.\\
We suppose the contrary, i.e., there exist $\delta>0$ such that
$H(n, m)\geq \delta$ for all sufficiently large integers $m$ and
$n$. Hence from the independence of the increments of $S$,
\[
\begin{split}
1&\geq \p\left(\bigcup_{m=n+1}^{\infty}C_{m}\right)\geq
\sum_{m=n+1}^{\infty}\p\left(C_{m}\bigcap\left(\bigcap_{j=n}^{m-1}C^{c}_{j}\right)\right)\\
&\geq\sum_{m=n+1}^{\infty}\p\Big(r^{m}>rh\Big(S_{r^{m}}\Big)\Big)H(n,m)
\geq \delta
 \sum_{m=n+1}^{\infty}\p\big(C_{m}\big),\\
\end{split}
\]
but since $\sum \p(C_{n})$ diverges, we see that our assertion is
true.\\
Now, we define
\[
\rho_{n_{l}, m_{l}}(x)\eqdef\p\Big(r^j\leq
rh\big(S_{r^j}-S_{r^{m_l-1}}+x\big)\textrm{ for, } n_{l}\leq j\leq
m_{l}-2\Big),\qquad x\geq 0,
 \]
and
\[
G(n_{l}, m_{l})\eqdef\p\Big(r^j\leq rh\big(S_{r^j}\big)\textrm{ for,
} n_{l}\leq j\leq m_{l}-1\Big).
\]
Since $h$ is increasing, we see that $\rho_{n_{l}, m_{l}}(x)$ is
increasing in $x$.\\
If we denote by $\mu$ and $\bar{\mu}$ the laws of $S_{1}$ and
$S_{1}-S_{r}$ respectively, by the scaling property we may express
$H(n_l,m_l)$ and $G(n_{l}, m_{l})$ as follows
\begin{align*}
H(n_l,m_l)&=\int_{0}^{+\infty}\bar{\mu}(\ud
x)\ind_{\big\{h(r^{m_{l}-1}x)\geq r^{m_{l}}\big\}}\rho_{n_{l},
m_{l}}(r^{m_{l}-1}x)\quad \textrm{and},\\
G(n_l,m_l)&=\int_{0}^{+\infty}\mu(\ud
x)\ind_{\big\{h(r^{m_{l}-1}x)\geq r^{m_{l}}\big\}}\rho_{n_{l},
m_{l}}(r^{m_{l}-1}x).
\end{align*}
In particular, we get that for $l$ sufficiently large
\[
H(n_l,m_l)\geq \rho_{n_{l}, m_{l}}(N)\int_{N}^{+\infty}\bar{\mu}(\ud
x)\qquad{for }\quad N\geq rC,
\]
where $C=\sup_{x\leq \beta}x/h(x)$. \\
Since $H(n_l,m_l)$ converges to $0$, as $l$ goes to $+\infty$ and
$\bar{\mu}$ does not depend on $l$, then $\rho_{n_{l},m_{l}}(N)$
also converges to $0$ when $l$ goes to $+\infty$, for every $N\geq
rC$.\\
On the other hand, we have
\[
G(n_{l}, m_{l})\leq \rho_{n_{l},m_{l}}(N)\int_{0}^{N}\mu(\ud x)+
\int_{N}^{\infty}\mu(\ud x),
\]
then letting $l$ and $N$ go to infinity, we get that  $G(n_{l},
m_{l})$ goes to $0$. Then, by (\ref{desteo}) we get that $\lim
\p(B_{n})=1$ and with this we finish the proof. \QED \noindent For
the integral tests at $+\infty$, we define $\mathcal{H}_{\infty}$,
the totality of positive increasing functions $h(t)$ on $(0,
\infty)$ that satisfy
\begin{itemize}
\item[i)] $\lim_{t\to\infty}h(t)=\infty$, and \item[ii)] there
exists $\beta>1$ such that
$\displaystyle\sup_{t>\beta}\displaystyle\frac{t}{h(t)}<\infty.$
\end{itemize}
The upper envelope of $X^{(x)}$ at $+\infty$ is given by the
following result.
\begin{theorem} Let $m\geq 0$ and $h\in\mathcal{H}_{\infty}$.
\begin{itemize}
\item[i)] If
\[
\int^{+\infty}F^\uparrow\left(\frac{t}{h(t)}\right)\frac{\ud
t}{t}<\infty,
\]
then for all $\epsilon > 0$ and for all $x\geq 0,$
\[
\p_{x}\Big(X_{t}>(1+\epsilon)h(t),\textrm{ i.o., as }t\to
+\infty\Big)=0.
\]
\item[ii)] Assume that (\ref{adhyp}) is satisfied. If
\[
\int^{+\infty}F^\uparrow\left(\frac{t}{h(t)}\right)\frac{\ud
t}{t}=\infty,
\]
then for all $\epsilon > 0$ and for all $x\geq 0$
\[
\p_{x}\Big(X_{t}>(1-\epsilon)h(t),\textrm{ i.o., as }t\to
+\infty\Big)=1.
\]
\end{itemize}
\end{theorem}
\textit{Proof:} We first consider the case where $x=0$. In this case
the proof of the tests at $+\infty$ is almost the same as that of
the tests at $0$. It is enough to apply the same arguments
to the sequence $x_{n}=r^{n}$, for $r>1$.\\
Now, we prove $(i)$ for any $x>0$. Let $h\in \mathcal{H}_{\infty}$
such that
$\int^{+\infty}F^\uparrow\left(\frac{t}{h(t)}\right)\frac{\ud t}{t}$
is finite. Let $x>0$ and $S_{x}=\inf\{t\geq
0\mbox{}:\mbox{}X^{(0)}_{t}\geq x\}$. Since clearly
\[
\int^{+\infty}F^\uparrow\left(\frac{t}{h(t-S_{x})}\right)\frac{\ud
t}{t}<\infty,
\]
from the Markov property at time $S_{x}$, we have for all
$\epsilon>0$
\begin{equation*}
\p_{0}\Big(X_{t}>(1+\epsilon)h(t-S_{x}), \textrm{ i. o., as } t\to
\infty\Big) =\p_{x}\Big(X_{t}>(1+\epsilon)h(t), \textrm{ i. o., as }
t\to \infty\Big)=0,
\end{equation*}
which proves part $(i)$.\\
Part $(ii)$ can be proved in the same way.\QED
\subsection{The general case.}
\begin{proposition}  Let $m>0$ and
$h\in\mathcal{H}_{0}$.
\begin{itemize}
\item[i)] If
\[
\int_{0^{+}}G\left(\frac{t}{h(t)}\right)\frac{\ud t}{t}<\infty,
\]
then for all $\epsilon > 0$
\[
\p\Big(X^{(0)}_{t}>(1+\epsilon)h(t),\textrm{ i.o., as } t\to
0\Big)=0.
\]
\item[ii)]If
\[
\int_{0^{+}}\overline{F}\left(\frac{t}{h(t)}\right)\frac{\ud t}{
t}=\infty,
\]
then for all $\epsilon > 0$
\[
\p\Big(X^{(0)}_t<(1-\epsilon)h(t),\textrm{ i.o., as }t\to 0\Big)=1.
\]
\end{itemize}
\end{proposition}
\textit{Proof:} Let $(x_{n})$ be a decreasing sequence which
converges to $0$. We define the events $A_{n}=\big\{\textrm{There
exists }t\in [S_{x_{n+1}}, S_{x_{n}})\textrm{ such that }
X^{(0)}_{t}>h(t)\big\}$. From the fact that $S_{x_{n}}$ tends to
$0$, a.s. when $n$ goes to $+\infty$, we see
\[
\Big\{X^{(0)}_{t}>h(t), \textrm{ i.o., as }t\to
0\Big\}=\limsup_{n\to +\infty} A_{n}.
\]
Since $h$ is an increasing function the following inclusion hold
\begin{equation}\label{inclusion}
A_{n}\subset \Big\{x_{n}>h\big(S_{x_{n+1}}\big)\Big\}.
\end{equation}
Now, we prove the convergent part. We choose $x_{n}=r^{n}$, for
$r<1$ and $h_{r}(t)=r^{-2}h(t)$. Since $h$ is increasing, we deduce
\[
\sum_{n} \p\Big(r^n>h_{r}\big(S_{r^{n+1}}\big)\Big)\leq
-\frac{1}{\log r}\int_{0}^{r}\p\Big(t>h\big(S_t\big)\Big)\frac{\ud
t}{t}.
\]
Replacing $h$ by $h_r$ in (\ref{inclusion}), we see that we can
obtain our result if
\[
 \int_{0}^{r}\p\Big(t>h\big(S_t\big)\Big)\frac{\ud t}{t}<\infty.
\]
From elementary calculations, we get
\[
\int_{0}^{r}\p\Big(t>h\big(S_t\big)\Big)\frac{\ud t}{t}=\e \left(
\int_{0}^{h^{-1}(r)}\ind_{\big\{t/r<S_1<t/h(t)\big\}}\frac{\ud
t}{t}\right),
\]
where $h^{-1}(s)=\inf\{t>0, h(t)>s\}$, the right inverse function of
$h$. Then, this integral converges if
\[
\int_{0}^{h^{-1}(r)}\p\left(S_1<\frac{t}{h(t)}\right)\frac{\ud t}{t}
<\infty.
\]
This proves part $(i)$.\\
The divergent part follows from the integral test for the upper
envelope of the future infimum of pssMp due to Pardo \cite{Pa}, (
see Theorem 1, part $(ii)$).\QED
\begin{proposition}  Let $m>0$ and
$h\in\mathcal{H}_{\infty}$.
\begin{itemize}
\item[i)] If
\[
\int^{+\infty }G\left(\frac{t}{h(t)}\right)\frac{\ud t}{t}<\infty,
\]
then for all $\epsilon > 0$ and for all $x\geq 0$
\[
\p\Big(X^{(x)}_{t}>(1+\epsilon)h(t),\textrm{ i.o., as } t\to +\infty
\Big)=0.
\]
\item[ii)]If
\[
\int^{+\infty}\overline{F}\left(\frac{t}{h(t)}\right)\frac{\ud t}{
t}=\infty,
\]
then for all $\epsilon > 0$ and for all $x\geq 0$
\[
\p\Big(X^{(x)}_t<(1-\epsilon)h(t),\textrm{ i.o., as }t\to
+\infty\Big)=1.
\]
\end{itemize}
\end{proposition}
\textit{Proof:} We first consider the case where $x=0$. In this case
the proof of the tests at $+\infty$ is almost the same as that of
the tests at $0$. It is enough to apply the same arguments to the
sequence $x_{n}=r^{n}$,
for $r>1$.\\
Now, we prove $(i)$ for any $x>0$. Let $h\in \mathcal{H}_{\infty}$
such that $\int^{+\infty}G\left(\frac{t}{h(t)}\right)\frac{\ud
t}{t}$ is finite. Let $x>0$ and $S_{x}$ and note by $\mu_{x}$ the
law of $X^{(0)}_{S_{x}}$. Since clearly
\[
\int^{+\infty}G\left(\frac{t}{h(t-S_{x})}\right)\frac{\ud
t}{t}<\infty,
\]
from the Markov property at time $S_{x}$, we have for all
$\epsilon>0$
\begin{equation}\label{lliMark}
\begin{split}
&\p_{0}\Big(X_{t}>(1+\epsilon)h(t-S_{x}), \textrm{ i. o., as } t\to
\infty\Big)\\
=&\int_{[x,+\infty)}\p_{y}\Big(X_{t}>(1+\epsilon)h(t), \textrm{ i.
o., as } t\to \infty\Big)\mu_{x}(\ud y)=0.
\end{split}
\end{equation}
If $x$ is an atom of $\mu_{x}$, then equality (\ref{lliMark}) shows
that
\[
\p\Big(X^{(x)}_{t}>(1+\epsilon)h(t), \textrm{ i. o., as } t\to
\infty\Big)=0
\]
and the result is proved. Suppose that $x$ is not an atom of
$\mu_{x}$. From Theorem 1 in \cite{CCh},  we know that
$X^{(0)}_{S_{x}}\ed xe^{\theta}$, where $\theta$ is a positive r.v.
such that
\[
\xi_{T_{z}}-z \xrightarrow[z\to +\infty]{(w)}\theta.
\]
Then from section 2 in \cite{CCh}, the law of $\theta$ is given by
\[
\p(\theta>t)=\e(\sigma_{1})\int_{(t, \infty)}s\nu(\ud s), \quad
t\geq 0,
\]
where $\sigma$ is the upward ladder height process associated with
$\xi$ and $\nu$ its L\'evy measure.\\
Hence, $\p(e^{\theta}>z)>0$ for $z>1$, and  for any $\alpha>0$,
$\mu_{x}(x, x+\alpha)>0$. Hence (\ref{lliMark}) shows that there
exists $y>x$ such that
\[
\p\Big(X^{(y)}_{t}>(1+\epsilon)h(t), \textrm{ i. o., as } t\to
\infty\Big)=0,
\]
for all $\epsilon>0$. The previous allows us to conclude part $(i)$.\\
Part $(ii)$ can be proved in the same way.\QED

\section{The regular case.}
In this section, we will assume that $m>0$. According to Chaumont
and Pardo \cite{CP} the law of $U_{1}$, the last passage time below
level $1$, is the same as $\nu I(\hat{\xi})$ where $\nu$ is a
positive random variable bounded above by $1$ and independent of the
exponential functional $I(\hat{\xi})$. Hence, we have the following
inequality
$\nu I(\hat{\xi})\leq I(\hat{\xi})$ a.s.\\
Now, let us define
\[
 \overline{F}_{\nu}(t):=\p\big(\nu I(\hat{\xi})<t\big),
\]
and suppose that
\begin{equation}\label{reg}
ct^{\beta}L(t)\leq \overline{F}(t)\leq \overline{F}_{\nu}(t)\leq
Ct^{\beta}L(t)\qquad\textrm{ as }\quad t\to 0,
\end{equation}
where $\beta>0$, $c$ and $C$ are two positive constants such that
$c\leq C$ and $L$ is a slowly varying function at $0$. An important
example included in this case is when
$\overline{F}$ and $\overline{F}_{\nu}$ are regularly varying functions at $0$. \\

\begin{proposition}
Under condition (\ref{reg}), we have that
\[
ct^{\beta}L(t)\leq G(t) \leq C_{\epsilon}t^{\beta}L(t)\qquad\textrm{
as }\quad t\to 0,
\]
where $C_{\epsilon}$ is a positive constant bigger than $C$.
\end{proposition}
\textit{Proof:} The lower bound is clear since $\overline{F}(t)\leq
G(t)$, for all $t\geq 0$ and our assumption. Now, let us define
$M^{(0)}_{t}=\sup_{0\leq s\leq t} X^{(0)}_{s}$ and fix $\epsilon
>0$. Then, by the Markov property and the fact that $J^{(x)}$ is an
increasing process, we have
\[
\begin{split}
\p_{0}\left(J_{1}> \frac{1-\epsilon}{t}\right)&\geq
\p_{0}\left(J_{1}>
\frac{1-\epsilon}{t}, M_{1}\geq \frac{1}{t}\right)\\
&=\e\left(S_{1/t}\leq 1, \p_{X^{(0)}_{S_{1/t}}}\left(J_{1-S_{1/t}}>
\frac{1-\epsilon}{t}\right)\right)\\
&\geq \e\left(S_{1/t}\leq 1, \p_{X^{(0)}_{S_{1/t}}}\left(J_{0}>
\frac{1-\epsilon}{t}\right)\right).\\
\end{split}
\]
Since $X^{(0)}_{S_{1/t}}\geq 1/t$ a.s., and the Lamperti
representation (\ref{lamp}), we deduce that
\[
\e\left(S_{1/t}\leq 1, \p_{X^{(0)}_{S_{1/t}}}\left(J_{0}>
\frac{1-\epsilon}{t}\right)\right) \geq \p\big(S_{1/t} <1
\big)\p\Big(\inf_{s\geq 0}\xi_{s}> \log (1-\epsilon)\Big).
\]
On the other hand, under the assumption that $\xi$ drifts towards
$+\infty$, we know from Section 2 of Chaumont and Doney \cite{CD}
(see also Proposition VI.17 in \cite{be}) that for all $\epsilon>0$
\[
K_{\epsilon}:=\p\Big(\inf_{s\geq 0}\xi_{s}> \log
(1-\epsilon)\Big)>0.
\]
Hence
\[
K^{-1}_{\epsilon}\p_{0}\left(J_{1}>
\frac{1-\epsilon}{t}\right)\geq\p\big(S_{1} <t \big)
\]
which implies that
\[
C K^{-1}_{\epsilon}\left(\frac{t}{1-\epsilon}\right)^{\beta}L(t)\geq
K^{-1}_{\epsilon}\p\left(U_{1}< \frac{t}{1-\epsilon}\right)\geq
\p\big(S_{1} <t \big) ,\qquad \textrm{as } \quad t\to 0,
\]
then the proposition is proved.\QED \noindent The next result give
us integral tests for the lower envelope of $S$ at $0$ and at
$\infty$, under condition (\ref{reg}).
\begin{theorem}
Under condition (\ref{reg}), the lower envelope of  $S$ at $0$ and
at $+\infty$ is as follows:
\begin{itemize}
\item[i)] Let $h\in\mathcal{H}^{-1}_{0}$, such that either $\lim_{x\to
0}h(x)/x=0$ or $\liminf_{x\to 0}h(x)/x>0$, then
\[
\p\Big(S_x<h(x), \textrm{ i.o., as } x\to 0\Big)= 0\textrm{ or } 1,
\]
according as
\[
\int_{0^+} \overline{F}\left(\frac{h(x)}{x}\right)\frac{\ud
x}{x}\qquad\textrm{is finite or infinite}.
\]\\
\item[ii)]Let $h\in\mathcal{H}^{-1}_{\infty}$, such that either
$\lim_{x\to +\infty}h(x)/x=0$ or $\liminf_{x\to +\infty}h(x)/x>0$,
then
\[
\p\Big(S_{x}<h(x), \textrm{ i.o., as } x\to \infty\Big)= 0\textrm{
or } 1,
\]
according as
\[
\int^{+\infty}\overline{F}\left(\frac{h(x)}{x}\right)\frac{\ud
x}{x}\qquad\textrm{is finite or infinite}.
\]
\end{itemize}
\end{theorem}
{\it Proof:} First let us check that under condition (\ref{reg}) we
have
\begin{equation}\label{equivre0}
\int_{0}^{\lambda}\overline{F}\left(\frac{h(x)}{x}\right)\frac{\ud
x}{x}<\infty\quad\textrm{ if and only if
}\quad\int_{0}^{\lambda}G\left(\frac{ h(x)}{x}\right)\frac{\ud
x}{x}<\infty.
\end{equation}
Since $\overline{F}(t)\leq G(t)$ for all $t\geq 0$, it is clear that
we only need to prove that
\[
\int_{0}^{\lambda}\overline{F}\left(\frac{h(x)}{x}\right)\frac{\ud
x}{x}<\infty\quad\textrm{ implies that
}\quad\int_{0}^{\lambda}G\left(\frac{ h(x)}{x}\right)\frac{\ud
x}{x}<\infty.
\]
From the hypothesis, either $\lim_{x\to 0}h(x)/x=0$ or
$\liminf_{x\to 0}h(x)/x>0$. In the first case, from condition
(\ref{reg}) there exists $\lambda>0$ such that, for every
$x<\lambda$
\[
c\left(\frac{h(x)}{x}\right)^{\beta}L\left(\frac{h(x)}{x}\right)\leq
\overline{F}\left(\frac{h(x)}{x}\right)\leq
C\left(\frac{h(x)}{x}\right)^{\beta}L\left(\frac{h(x)}{x}\right).
\]
Since, we suppose that
$\int_{0}^{\lambda}\overline{F}\left(\frac{h(x)}{x}\right)\frac{\ud
x}{x}$ is finite, then
\[
\int_{0}^{\lambda}\left(\frac{h(x)}{x}\right)^{\beta}L\left(\frac{h(x)}{x}\right)\frac{\ud
x}{x}<\infty, \] hence from Proposition 8, we get that
$\int_{0}^{\lambda}G\left(\frac{h(x)}{x}\right)\frac{\ud x}{x}$ is
also finite. In the second case, since for any $0<\delta<\infty$,
$\p\big(I<\delta\big)>0$, and $\liminf_{x\to 0}h(x)/x>0$, we have
for any $y$
 \begin{equation}\label{inprobreg}
0<\p\left(I<\liminf_{x\to 0}\frac{h(x)}{x}\right)<\p\left(
I<\frac{h(y)}{y}\right).
\end{equation}
Hence, since for every $t\geq 0$, $\overline{F}(t)\leq G(t)$, we
deduce that
\[
\int_{0}^{\lambda}\overline{F}\left(\frac{h(x)}{x}\right)\frac{\ud
x}{x}= \int_{0}^{\lambda}G\left(\frac{h(x)}{x}\right)\frac{\ud
x}{x}=\infty.
\]
Now, let us  check that for any constant $\beta>0$,
\begin{equation}\label{equivre}
\int_{0}^{\lambda}\overline{F}\left(\frac{h(x)}{x}\right)\frac{\ud
x}{x}<\infty\quad\textrm{ if and only if
}\quad\int_{0}^{\lambda}\overline{F}\left(\frac{\beta
h(x)}{x}\right)\frac{\ud x}{x}<\infty,
\end{equation}
Again, from the hypothesis either $\lim_{x\to 0}h(x)/x=0$ or
$\liminf_{x\to 0}h(x)/x>0$. In the first case, we deduce
(\ref{equivre}) from (\ref{reg}). In the
second case, from (\ref{inprobreg}) both of the integrals in (\ref{equivre}) are infinite.\\
Next, it follows from Proposition 4 part $(i)$ and (\ref{equivre0})
that if $\int_{0^+}\overline{F}\left(\frac{h(x)}{x}\right)\frac{\ud
x}{x}$ is finite, then for all $\epsilon>0$,
\[
\p\big(S_x<(1-\epsilon)h(x), \textrm{ i.o., as }x\to 0\big)=0.
\]
If $\int_{0^+} \overline{F}\left(\frac{h(x)}{x}\right)\frac{\ud
x}{x}$ diverges, then from Proposition 4 part $(ii)$ that for all
$\epsilon>0$,
\[
\p\big(S_x<(1+\epsilon)h(x), \textrm{ i.o., as }x\to 0\big)=1.
\]
Then $(\ref{equivre})$  allows
us to drop $\epsilon$ in this implications.\\
The tests at $+\infty$ are proven through the same way.\QED
\noindent From the previous Theorem, we deduce that under condition
(\ref{reg}) the first and the last passage time processes have the
same lower functions ( see Theorem 5 in \cite{Pa}).

\begin{theorem}
Under condition (\ref{reg}), the upper envelope of the pssMp  at $0$
and at $+\infty$ is as follows:
\begin{itemize}
\item[i)] Let $h\in\mathcal{H}_{0}$, such that either $\lim_{t\to
0}t/h(t)=0$ or $\liminf_{t\to 0}t/h(t)>0$, then
\[
\p\Big(X^{(0)}_{t}>h(t), \textrm{ i.o., as } t\to 0\Big)= 0\textrm{
or } 1,
\]
according as
\[
\int_{0^+} \overline{F}\left(\frac{t}{h(t)}\right)\frac{\ud
t}{t}\qquad\textrm{is finite or infinite}.
\]\\
\item[ii)]Let $h\in\mathcal{H}_{\infty}$, such that either
$\lim_{t\to +\infty}t/h(t)=0$ or $\liminf_{t\to +\infty}t/h(t)>0$,
then for all $x\geq 0$
\[
\p\Big(X^{(x)}_{t}>h(t), \textrm{ i.o., as } t\to \infty\Big)=
0\textrm{ or } 1,
\]
according as
\[
\int^{+\infty}\overline{F}\left(\frac{t}{h(t)}\right)\frac{\ud
t}{t}\qquad\textrm{is finite or infinite}.
\]
\end{itemize}
\end{theorem}
{\it Proof:} We prove this result by following the same arguments as
the proof of the previous Theorem.\QED

Note that from this result, we deduce that under condition
(\ref{reg}) a pssMp and its future infimum have the same upper
functions ( see Theorem 6 in \cite{Pa}).
\section{The log-regular case.}
In this section, we also assume that $m>0$. Here, we will study to
types of behaviour of $\overline{F}$ and $ \overline{F}_{\nu}$, both
types of behaviour allow us to obtain laws of the iterated logarithm
for the upper envelope of $X^{(0)}$. The first type of behaviour
that we will consider is when  $\log \overline{F}$ and $\log
\overline{F}_{\nu}$ are regularly varying at $0$, i.e.
\begin{equation}\label{logreg}
-\log  \overline{F}_{\nu}(1/t)\sim-\log  \overline{F}(1/t)\sim
\lambda t^{\delta}L(t),\quad \textrm{ as } t\to +\infty,
\end{equation}
where $\lambda>0$, $\delta>0$ and $L$ is a slowly varying function
at
$+\infty$.\\
The second type of behaviour that we will consider is when  $\log
\overline{F}$ and $\log \overline{F}_{\nu}$ satisfy that
\begin{equation}\label{loglogreg}
-\log \overline{F}_{\nu}(1/t)\sim-\log \overline{F}(1/t)\sim  K
(\log t)^{\gamma},\quad \textrm{ as } t\to +\infty,
\end{equation}
where $K$ and $\gamma$ are strictly positive constants.\\
\subsection{Laws of the iterated logarithm for pssMp.}
\begin{proposition}
Under condition (\ref{logreg}), the tail probability of $S_1$
satisfies
\begin{equation}\label{suplogreg}
-\log G(1/t) \sim \lambda t^{\delta}L(t)\qquad\textrm{ as }\quad
t\to +\infty.
\end{equation}
Similarly, under condition (\ref{loglogreg}), the tail probability
of $S_1$ satisfies
\begin{equation}\label{suploglogreg}
-\log G(1/t) \sim  K( \log t)^{\gamma}\qquad\textrm{ as }\quad t\to
+\infty.
\end{equation}
\end{proposition}
\textit{Proof:} First, we prove the upper bound of
(\ref{suplogreg}). With the same notation as in the proof of
Proposition 8, we see
\[
-\log \p\Big(\nu I(\hat{\xi})<1/t\Big)=-\log
\p_{0}\Big(J_{1}>t\Big)\geq -\log \p_{0}\Big(M_{1}>t\Big),
\]
which implies
\[
1\geq \limsup_{t\to \infty}\frac{-\log
\p_{0}\Big(M_{1}>t\Big)}{\lambda t^{\delta}L(t)},
\]
and since $\p_{0}\big(M_{1}>t\big)=\p\big(S_1<1/t\big)$, we
get the upper bound.\\
Now, fix $\epsilon >0$. From the proof of Proposition 8, we have
\[
\p_{0}\big(J_{1}> (1-\epsilon)t\big)\geq \p\big(S_{t} <1
\big)\p\Big(\inf_{s\geq 0}\xi_{s}> \log(1-\epsilon)\Big).
\]
On the other hand, we know
\[
K_{\epsilon}:=\p\Big(\inf_{s\geq 0}\xi_{s}> \log(1-\epsilon)\Big)>0,
\]
Hence,
\[
-\log \p_{0}\big(J_{1}> (1-\epsilon)t\big)\leq -\log\p\Big(S_{1}
<1/t\Big)
 -\log K_{\epsilon},
\]
which implies the following lower bound
\[
(1-\epsilon)^{\delta}\leq \liminf_{t\to \infty}\frac{-\log
\p\Big(S_{1} <1/t\Big)}{\lambda t^{\delta}L(t)},
\]
and since $\epsilon$ can be chosen arbitrarily small,
(\ref{suplogreg}) is proved.\\
The upper bound of tail behaviour (\ref{suploglogreg}) is proven
through the same way. For the lower bound, we follow the same
arguments as above and we get that
\[
-\log \p_{0}\big(J_{1}> (1-\epsilon)t\big)\leq -\log\p\Big(S_{1}
<1/t\Big)
 -\log K_{\epsilon},
\]
which implies
\[ 1=\liminf_{t\to\infty} \left(\frac{\log (1-\epsilon)t}{\log t}\right)^{\gamma}\leq \liminf_{t\to
\infty}\frac{-\log \p\Big(S_{1} <1/t\Big)}{K (\log t)^{\gamma}},
\]
then the proposition is proved. \QED \noindent The following result
gives us laws of the iterated logarithm for the first passage time
process when condition (\ref{logreg}) is satisfied.\\
Define the functions
\begin{equation*}
\varphi(x):= \frac{x}{\inf\big\{s:1/ \overline{F}(1/s)>|\log
x|\big\}},\quad x>0, \quad x\neq 1,
\end{equation*}
and
\begin{equation*}
\vartheta(t):= \frac{t^2}{\varphi(t)},\quad t>0, \quad t\neq 1.
\end{equation*}
\begin{theorem}Under condition (\ref{logreg}), we have the following law of the iterated logarithm for $S$:
\[
\limsup_{x\to
0}\frac{S_{x}}{\varphi(x)}=1\quad\textrm{and}\quad\limsup_{x\to
\infty}\frac{S_{x}}{\varphi(x)}=1\quad \textrm{almost surely.}
\]
The upper envelope of pssMp, under condition (\ref{logreg}), are
described by the following law of the iterated logarithm:
\begin{itemize}
\item[i)]\[ \limsup_{t\to 0}\frac{X^{(0)}_{t}}{\vartheta(t)}=1,\qquad
\textrm{almost surely.}
\]
\item[ii)]For all $x\geq 0$,
\[ \limsup_{t\to
+\infty}\frac{X^{(x)}_{t}}{\vartheta(t)}=1,\qquad \textrm{almost
surely.}
\]
\end{itemize}
\end{theorem}
{\it Proof:} This Theorem is a consequence of Propositions 4, 5,6,7
and 9, and it is proven in the same way as Theorem 4 in \cite{CP},
we only need to emphasize that we can replace $\log G$ by $\log
\overline{F}$, since they are asymptotically equivalent. \QED

Note that under condition (\ref{logreg}) a pssMp and its future
infimum satisfy the same law of the iterated logarithm (see Theorem
8 in \cite{Pa}) but they do
not necessarily have the same upper functions.\\
Similarly, under condition (\ref{loglogreg}) we may establish laws
of the iterated logarithm for the upper envelope of pssMp and their
future infimum. In this direction, let us define
\[
\phi(x):=x\exp\Big\{-\big(K^{-1}\log |\log
x|\big)^{1/\gamma}\Big\},\quad x>0,\quad x\neq 1,
\]
and
\[
\Phi(t):=\frac{t^2}{\phi(t)},\quad t>0,\quad t\neq 1.
\]
We recall that $J^{(x)}=(J^{(x)}_t, t\ge 0)$ is the future infimum
process of $X^{(x)}$, for $x\geq 0$, where $J^{(x)}_{t}=\inf_{s\geq
t} X^{(x)}$.
\begin{theorem}Under condition (\ref{loglogreg}),we have the following
laws of the iterated logarithm:
\begin{itemize}
\item[i)] For the first passage time, we have
\[ \limsup_{x\to 0}\frac{S_{x}}{\phi(x)}=1, \qquad
\limsup_{x\to \infty}\frac{S_{x}}{\phi(x)}=1\quad \textrm{almost
surely.}\]
\item[ii)] For the last passage time, we have
\[
\limsup_{x\to 0}\frac{U_{x}}{\phi(x)}=1 \quad\textrm{and}\quad
\limsup_{x\to +\infty}\frac{U_{x}}{\phi(x)}=1\quad \textrm{almost
surely.}
\]
\end{itemize}
The upper envelope of pssMp and their future infimum processes,
under condition (\ref{loglogreg}), are described by the following
laws of the iterated logarithm:
\begin{itemize}
\item[iii)]\[ \limsup_{t\to 0}\frac{X^{(0)}_{t}}{\Phi(t)}=1
\quad\textrm{and}\quad \limsup_{t\to
0}\frac{J^{(0)}_{t}}{\Phi(t)}=1\quad \textrm{almost surely.}
\]
\item[iv)]For all $x\geq 0$,
\[ \limsup_{t\to
+\infty}\frac{X^{(x)}_{t}}{\Phi(t)}=1 \quad\textrm{and}\quad
\limsup_{t\to 0}\frac{J^{(x)}_{t}}{\Phi(t)}=1\quad \textrm{almost
surely.}
\]
\end{itemize}
\end{theorem}
{\it Proof}: We first prove part $(i)$ for small times. Note that it
is easy to check that both $\phi(x)$ and $\phi(x)/x$ are increasing
in a neighbourhood of 0, moreover the function $\phi(x)/x$ is
bounded
by 1, for $x\in [0, 1)$.\\
From condition (\ref{loglogreg}) and Proposition 9, we have for all
$k_1<1$ and $k_2>1$ and for all $t$ sufficiently large,
\begin{eqnarray*}
k_1 K (\log t)^{\gamma}\le -\log G(1/t)\le k_2K (\log t)^{\gamma},
\end{eqnarray*}
so that for $\phi$ defined above,
\[
k_1 \log |\log t|\le -\log G\left(\phi(x)/x\right)\le k_2\log |\log
t|,
\]
hence
\[
G\left(\frac{\phi(x)}{x}\right)\ge  (|\log t|)^{-k_2}\,.
\]
Since $k_2>1$, we obtain the convergence of the integral
\[
\int_{0+}G\left(\frac{\phi(x)}{x}\right)\,\frac{\ud t}{t}\,,
\]
which proves that for all $\varepsilon>0$,
\[
\p(S_x<(1-\varepsilon)\phi(x),\mbox{i.o., as $x\rightarrow0$})=0
\] from Proposition 4 part $(i)$. The
divergent part is proven through the same way so that from
Proposition 4 part $(ii)$, one has for all $\varepsilon>0$,
\[\p(S_x<(1+\varepsilon)\phi(x),\mbox{i.o., as
$x\rightarrow0$})=1
\] and the conclusion follows.\\
Condition (\ref{loglogreg}) implies that $\phi(x)$ is increasing in
a neighbourhood of $+\infty$ whereas $\phi(x)/x$ is decreasing in a
neighbourhood of $+\infty$. Hence, the proof of the result at
$+\infty$ is done through the same way as at 0, by using Proposition
5.\\
The parts (ii), (iii) and (iv) can be proved following the same
arguments, it is enough to specify that the laws of the iterated
logarithm of the last passage process and the future infimum process
will use the integral tests found by Pardo \cite{Pa}, (see Theorems
1,2,3 and 4 in \cite{Pa})\QED
\subsection{The case with no positive jumps.}
Here, we suppose that the pssMp $X^{(x)}$ has no positive jumps. Our
next result is a remarkable asymptotic property of this type of
pssMp. This Theorem means in particular that if there exists a
positive function that describes the upper envelope of $X^{(x)}$ by
a law of the iterated logarithm then the same function describes the
upper envelope of the future infimum of $X^{(x)}$ and the pssMp
$X^{(x)}$ reflected at its future infimum.
\begin{theorem} Let us suppose that
\[
\limsup_{t\to 0}\frac{X^{(0)}_t}{\Lambda(t)}=1\qquad\textrm{ almost
surely,}
\]
where $\Lambda$ is a positive function such that $\Lambda(0)=0$,
then
\[
\limsup_{t\to 0}\frac{J^{(0)}_t}{\Lambda(t)}=1\quad \textrm{ and
}\quad \limsup_{t\to
0}\frac{X^{(0)}_t-J^{(0)}_t}{\Lambda(t)}=1\quad\textrm{ almost
surely.}
\]
Moreover, if for all $x\geq 0$
\[
\limsup_{t\to +\infty}\frac{X^{(x)}_t}{\Lambda(t)}=1\qquad\textrm{
almost surely,}
\]
where $\Lambda$ is a positive function such that
$\lim_{t\to+\infty}\Lambda(t)=+\infty$, then
\[
\limsup_{t\to +\infty}\frac{J^{(x)}_t}{\Lambda(t)}=1\quad \textrm{
and }\quad \limsup_{t\to
+\infty}\frac{X^{(0)}_t-J^{(0)}_t}{\Lambda(t)}=1\quad\textrm{ almost
surely.}
\]
\end{theorem}
\textit{Proof: } First, we prove the result for large times. Let
$x\geq 0$. Since $J^{(x)}_{t}\leq X^{(x)}_{t}$ for every $t\geq 0$
and our hypothesis, then it is clear
\[
\limsup_{t\to +\infty}\frac{J^{(x)}_t}{\Lambda(t)}\leq 1\qquad
\textrm{ almost surely.}
\]
Now, fix $\epsilon\in(0,1/2)$ and define
\[
R_n =\inf\left\{s\geq n: \frac{X^{(x)}_{s}}{\Lambda(s)}\geq
(1-\epsilon)\right\}.
\]
From the above definition, it is clear that $R_{n}\geq n$ and that
$R_{n}$ diverge a.s. as $n$ goes to $+\infty$. From our hypothesis,
we
deduce that $R_{n}$ is finite, a.s.\\
Now, since $X^{(x)}$ has no positive jumps and applying the strong
Markov property and Lamperti representation (\ref{lamp}), we have
\[
\begin{split}
\p\left(\frac{J^{(x)}_{R_n}}{\Lambda(R_n)}\geq
(1-2\epsilon)\right)&=\p\left(J^{(x)}_{R_n}\geq
\frac{(1-2\epsilon)X^{(x)}_{R_{n}}}{(1-\epsilon)}\right)\\
&=\e\left( \p\left(J^{(x)}_{R_{n}}\geq
\frac{(1-2\epsilon)X^{(x)}_{R_{n}}}{(1-\epsilon)}\bigg|X^{(x)}_{R_{n}} \right)\right)\\
&=\p\left(\inf_{t\geq 0}\xi\geq
\log\frac{(1-2\epsilon)}{(1-\epsilon)}\right)=cW\left(\log
\frac{1-\epsilon}{1-2\epsilon}\right)>0,
\end{split}
\]
where $W:[0,+\infty)\to[0,+\infty)$ is the unique absolutely
continuous increasing function with Laplace exponent
\[
\int_{0}^{+\infty} e^{-\lambda x} W(x)\ud x =\frac{1}{\psi(\lambda)}
\quad \textrm{ for }\lambda >0,
\]
and $c=1/W(+\infty)$, (see Bertoin \cite{be} Theorem VII.8).\\
Since $R_{n}\geq n$,
\[
\p\left(\frac{J^{(x)}_{t}}{\Lambda(t)}\geq (1-2\epsilon), \textrm{
for some } t\geq
n\right)\geq\p\left(\frac{J^{(x)}_{R_n}}{\Lambda(R_n)}\geq
(1-2\epsilon)\right).
\]
Therefore, for all $\epsilon\in (0, 1/2)$,
\[
\p\left(\frac{J^{(x)}_{t}}{\Lambda(t)}\geq (1-2\epsilon),
\textrm{i.o., as }t\to +\infty\right)\geq\lim_{n\to
+\infty}\p\left(\frac{J^{(x)}_{R_n}}{\Lambda(R_n)}\geq
(1-2\epsilon)\right)>0.
\]
The event of the left hand side is in the upper-tail sigma-field
$\cap_{t}\sigma\{X^{(x)}_{s}:s\geq t\}$ which is trivial, then
\[
 \limsup_{t\to +\infty}\frac{J^{(x)}_t}{\Lambda(t)}\geq
1-2\epsilon\qquad \textrm{ almost surely.}
\]
The proof of part $(ii)$ is very similar, in fact
\[
\begin{split}
\p\left(\frac{X^{(x)}_{R_n}-J^{(x)}_{R_n}}{\Lambda(R_n)}\geq
(1-2\epsilon)\right)&=\p\left(J^{(x)}_{R_n}\leq
\frac{\epsilon X^{(x)}_{R_{n}}}{1-\epsilon}\right)\\
&=\e\left( \p\left(J^{(x)}_{R_{n}}\leq
\frac{\epsilon X^{(x)}_{R_{n}}}{1-\epsilon}\bigg| X^{(x)}_{R_{n}}\right)\right)\\
&=\p\left(\inf_{t\geq 0}\xi\leq
\log\frac{\epsilon}{1-\epsilon}\right)=1-cW\left(\log
\frac{1-\epsilon}{\epsilon}\right)>0.
\end{split}
\]
Since $R_{n}\geq n$,
\[
\p\left(\frac{X^{(x)}_t-J^{(x)}_{t}}{\Lambda(t)}\geq (1-2\epsilon),
\textrm{ for some } t\geq
n\right)\geq\p\left(\frac{X^{(x)}_{R_n}-J^{(x)}_{R_n}}{\Lambda(R_n)}\geq
(1-2\epsilon)\right).
\]
Therefore, for all $\epsilon\in (0, 1/2)$,
\[
\p\left(\frac{X^{(x)}_t-J^{(x)}_{t}}{\Lambda(t)}\geq (1-2\epsilon),
\textrm{i.o., as }t\to \infty\right)\geq \lim_{n\to
+\infty}\p\left(\frac{X^{(x)}_{R_n}-J^{(x)}_{R_n}}{\Lambda(R_n)}\geq
(1-2\epsilon)\right)>0.
\]
The event of the left hand side of the above inequality is in the
upper-tail sigma-field $\cap_{t}\sigma\{X^{(x)}_{s}:s\geq t\}$ which
is trivial and this establishes part (ii) for large times.\\
In order to prove the LIL for small times, we now define the
following stopping time
\[
R_{n}=\inf\left\{\frac{1}{n}\leq s :
\frac{X^{(0)}_{s}}{\Lambda(s)}\geq (1-\epsilon)\right\}.
\]
Following same arguments as above, we get that for a fixed
$\epsilon\in(0,1/2)$ and $n$ sufficiently large
\[
\p\left(\frac{J^{(0)}_{R_n}}{\Lambda(R_n)}\geq
(1-2\epsilon)\right)>0\quad\textrm{ and
}\quad\p\left(\frac{X^{(0)}_{R_n}-J^{(0)}_{R_n}}{\Lambda(R_n)}\geq
(1-2\epsilon)\right)>0.
\]
Next, we note
\[
\p\left(\frac{J^{(0)}_{R_p}}{\Lambda(R_p)}\geq (1-2\epsilon),
\textrm{ for some } p\geq
n\right)\geq\p\left(\frac{J^{(0)}_{R_n}}{\Lambda(R_n)}\geq
(1-2\epsilon)\right),
\]
and
\[
\p\left(\frac{X^{(0)}_{R_p}-J^{(0)}_{R_p}}{\Lambda(R_p)}\geq
(1-2\epsilon), \textrm{ for some } p\geq
n\right)\geq\p\left(\frac{X^{(0)}_{R_n}-J^{(0)}_{R_n}}{\Lambda(R_n)}\geq
(1-2\epsilon)\right).
\]
Since $R_n$ converge a.s. to $0$ as $n$ goes to $\infty$, the
conclusion follows taking the limit when $n$ goes towards to
$+\infty$\QED
 \noindent
Hence, when $F^\uparrow$ satisfies condition (\ref{logreg}) we have
the following laws of the iterated logarithm for the future infimum
of $X^{(x)}$ and the pssMp $X^{(x)}$ reflected at its future
infimum.
\begin{corollary}
Under condition (\ref{logreg}), we have the following laws of the
iterated logarithm:
\begin{itemize}
\item[i)] for all $x\geq 0$
\[ \limsup_{t\to
0}\frac{J^{(0)}_t}{\vartheta(t)}=1\quad \textrm{ and }\quad
\limsup_{t\to +\infty}\frac{J^{(x)}_t}{\vartheta(t)}=1 \quad\textrm{
almost surely,}
\]\\
\item[ii)] for all $x\geq 0$
\[
\limsup_{t\to 0}\frac{X^{(0)}_t-J^{(0)}_t}{\vartheta(t)}=1\quad
\textrm{ and }\quad \limsup_{t\to
+\infty}\frac{X^{(x)}_t-J^{(x)}_t}{\vartheta(t)}=1\quad\textrm{
almost surely,}
\]
\end{itemize}
\end{corollary}

\section{Bessel processes}
Recall that Bessel processes are the only continuous positive
self-similar Markov processes. Recall also that a Bessel process of
dimension $\delta\geq 0$ with starting point $x\geq 0$ is the
diffusion $R$ whose square satisfies the stochastic differential
equation
\begin{equation}\label{carredube}
R^{2}_{t}=x^{2}+2\int_{0}^{t} R_{s}\ud \beta_{s}+\delta t, \qquad
t\geq 0,
\end{equation}
where $\beta$ is a standard Brownian Motion.\\
Now, we define $\xi=(2(B_{t}+at), t\geq 0)$, where $B$ is a standard
Brownian motion and $a\geq 0$. By the Lamperti representation, we
know that we can define a pssMp starting from $x>0$, such that
\[
X^{(x)}_{t}=x\exp\{\xi_{\tau(t/x)}\}, \quad t\ge 0.
\]
Applying the It\^o's formula and Dubins-Schwartz's Theorem (see for
instance Revuz and Yor \cite{RY}), we see that $X^{(x)}$ satisfy
(\ref{carredube}) with $\delta=2(a+1)$. Obviously, $\xi$ satisfies
the conditions under which we can define $X^{(x)}$ when $x=0$. When
$a>0$, it is clear that $X^{(x)}$ is transient and when $a=0$, it is
also clear that the process $X^{(x)}$ oscillates.\\
Gruet and Shi \cite{GS} proved that there exist a finite constant
$K>1$, such that for any $0<s\leq 2,$
\begin{equation}\label{desibes}
K^{-1}s^{1-\delta/2}\exp\left\{-\frac{1}{2s}\right\}\leq
\p(S_{1}<s)\leq K s^{1-\delta/2}\exp\left\{-\frac{1}{2s}\right\}.
\end{equation}
Hence we establish the following integral test for the lower
envelope of the first passage time process of the squared Bessel
process $X^{(0)}$.
\begin{theorem}
Let $h\in\mathcal{H}^{-1}_{0}$ and $\delta\ge 2$,
\begin{itemize}
\item[i)] If
\[
\int_{0^{+}}\left(\frac{t}{h(t)}\right)^{\frac{\delta-2}{2}}\exp\left\{-\frac{t}{2h(t)}\right\}\frac{\ud
t}{t}<\infty,
\]
then for all $\epsilon > 0$
\[
\p\Big(S_t<(1-\epsilon)h(t),\textrm{ i.o., as }t\to 0\Big)=0.
\]
\item[ii)] If
\[
\int_{0^{+}}\frac{\ud
t}{t}\left(\frac{t}{h(t)}\right)^{\frac{\delta-2}{2}}\exp\left\{-\frac{t}{2h(t)}\right\}=\infty,
\]
then for all $\epsilon > 0$
\[
\p\Big(S_t<(1+\epsilon)h(t),\textrm{ i.o., as }t\to 0\Big)=1.
\]
\end{itemize}
\end{theorem}
\textit{Proof: } The proof of this Theorem is a simple application
of (\ref{desibes}) to Proposition 2.\QED Similarly, we have the same
integral test for large times.
\begin{theorem}
Let $h\in\mathcal{H}^{-1}_{\infty}$ and $\delta\ge 2$,
\begin{itemize}
\item[i)] If
\[
\int^{+\infty}\left(\frac{t}{h(t)}\right)^{\frac{\delta-2}{2}}\exp\left\{-\frac{t}{2h(t)}\right\}\frac{\ud
t}{t}<\infty,
\]
then for all $\epsilon > 0$
\[
\p\Big(S_t<(1-\epsilon)h(t),\textrm{ i.o., as }t\to +\infty\Big)=0.
\]
\item[ii)] If
\[
\int_{+\infty}\frac{\ud
t}{t}\left(\frac{t}{h(t)}\right)^{\frac{\delta-2}{2}}\exp\left\{-\frac{t}{2h(t)}\right\}=\infty,
\]
then for all $\epsilon > 0$
\[
\p\Big(S_t<(1+\epsilon)h(t),\textrm{ i.o., as }t\to +\infty\Big)=1.
\]
\end{itemize}
\end{theorem}
From these integral tests, we get the following law of the iterated
logarithm
\[
\liminf_{t\to 0} S_{t}\frac{2\log |\log t|}{t}=1
\quad\textrm{and}\quad \liminf_{t\to +\infty} S_{t}\frac{2\log \log
t}{t}=1 \quad \textrm{almost surely.}
\]
For the upper envelope of $X^{(0)}$, we have the following
integral tests.
\begin{theorem}
Let $h\in\mathcal{H}_{0}$ and $\delta\ge 2$,
\begin{itemize}
\item[i)] If
\[
\int_{0^{+}}\left(\frac{h(t)}{t}\right)^{\frac{\delta-2}{2}}\exp\left\{-\frac{h(t)}{2t}\right\}\frac{\ud
t}{t}<\infty,
\]
then for all $\epsilon > 0$
\[
\p\Big(X^{(0)}_t>(1+\epsilon)h(t),\textrm{ i.o., as }t\to 0\Big)=0.
\]
\item[ii)] If
\[
\int_{0^{+}}\left(\frac{h(t)}{t}\right)^{\frac{\delta-2}{2}}\exp\left\{-\frac{h(t)}{2t}\right\}\frac{\ud
t}{t}=\infty,
\]
then for all $\epsilon > 0$
\[
\p\Big(X^{(0)}_t>(1-\epsilon)h(t),\textrm{ i.o., as }t\to 0\Big)=1.
\]
\end{itemize}
\end{theorem}
\textit{Proof: } The proof of this Theorem follows from a simple
application of (\ref{desibes}) to Theorem 2. The proof of the
additional hypothesis (\ref{adhyp}), is clear from
(\ref{desibes}).\QED Similarly, we have the same integral tests for
large times.
\begin{theorem}
Let $h\in\mathcal{H}_{\infty}$ and $\delta\ge 2$,
\begin{itemize}
\item[i)] If
\[
\int^{+\infty}\left(\frac{h(t)}{t}\right)^{\frac{\delta-2}{2}}\exp\left\{-\frac{h(t)}{2t}\right\}\frac{\ud
t}{t}<\infty,
\]
then for all $\epsilon > 0$ and for all $x\geq 0$
\[
\p\Big(X^{(x)}_t>(1+\epsilon)h(t),\textrm{ i.o., as }t\to
+\infty\Big)=0.
\]
\item[ii)] If
\[
\int^{+\infty}\left(\frac{h(t)}{t}\right)^{\frac{\delta-2}{2}}\exp\left\{-\frac{h(t)}{2t}\right\}\frac{\ud
t}{t}=\infty,
\]
then for all $\epsilon > 0$ and for all $x\geq 0$
\[
\p\Big(X^{(x)}_t>(1-\epsilon)h(t),\textrm{ i.o., as }t\to
+\infty\Big)=1.
\]
\end{itemize}
\end{theorem}
Recall from the Kolmogorov and Dvoretzky-Erd\"os (KDE for short)
integral test (Theorem 1) that for  $h$  a nondecreasing, positive
and unbounded function as $t$ goes to $+\infty$, the upper envelope
of $X^{(0)}$ at $0$ may be described as follows:
\[
\p\big(X^{(0)}_{t}> h(t), \textrm{ i.o., as } t\to 0\big)=0 \textrm{
or }1,
\]
according as,
\[
\int_{0}\left(\frac{h(t)}{t}\right)^{\frac{\delta}{2}}\exp\left\{-\frac{h(t)}{2t}\right\}\frac{\ud
t}{t} \qquad\textrm{is finite or infinite.}
\]
Note that the class of functions that satisfy the divergent part of
Theorems 11 and 12 implies the divergent part of the KDE integral
test, hence $\epsilon$ can also take the value $0$. The convergent
part of the KDE integral test obviously implies the convergent part
of Theorems 11 and 12.\\
From these integral tests, we get the following law of the iterated
logarithm
\[
\limsup_{t\to 0} \frac{X_{t}}{2t\log |\log t|}=1
\quad\textrm{and}\quad \limsup_{t\to +\infty} \frac{X^{(x)}}{2t\log
\log t}=1 \quad \textrm{almost surely,}
\]
for $x\geq 0$.\\

\section{Examples.}
{\bf Example 1.} The first example that we will consider here, is
the stable subordinator. Let $X^{(0)}$ be the stable subordinator
with index $\alpha\in(0,1)$. From Zolotarev \cite{zo1}, we know that
there exists $k$ a positive constant such that
\[
\p_0(X_1>x)\sim kx^{-\alpha} \qquad x\to +\infty,
\]
and since a subordinator is an increasing process, then
\[
\p(S_1<x)\sim kx^{\alpha} \qquad x\to 0,
\]
From Example 3 in Section 7.1 in \cite{Pa}, we also have that
\[
\p\Big(I(\hat{\xi})<x\Big)\sim mk\alpha x^{\alpha+1} \qquad x\to 0.
\]
Hence from Theorems $4$ and $5$, we obtain the following
corollaries.
\begin{corollary}
The lower envelope of  $S$, the first passage time of the stable
subordinator $X^{(0)}$ with index $\alpha\in(0,1)$ at $0$ and at
$+\infty$ is as follows:
\begin{itemize}
\item[i)] Let $h\in\mathcal{H}^{-1}_{0}$, such that either $\lim_{x\to
0}h(x)/x=0$ or $\liminf_{x\to 0}h(x)/x>0$, then
\[
\p\Big(S_x<h(t), \textrm{ i.o., as } x\to 0\Big)= 0\textrm{ or } 1,
\]
according as
\[
\int_{0^+} \left(\frac{h(x)}{x}\right)^{\alpha}\frac{\ud
x}{x}\quad\textrm{is finite } \quad\textrm{or}\quad \int_{0^+}
\left(\frac{h(x)}{x}\right)^{\alpha+1}\frac{\ud
x}{x}\qquad\textrm{is infinite}.
\]
\item[ii)]Let $h\in\mathcal{H}^{-1}_{\infty}$, such that either
$\lim_{x\to +\infty}h(x)/x=0$ or $\liminf_{x\to +\infty}h(x)/x>0$,
then
\[
\p\Big(S_{x}<h(x), \textrm{ i.o., as } x\to \infty\Big)= 0\textrm{
or } 1,
\]
according as
\[
\int^{+\infty}\left(\frac{h(x)}{x}\right)^\alpha\frac{\ud
x}{x}\qquad\textrm{is finite }\quad \textrm{or}\quad
\int^{+\infty}\left(\frac{h(x)}{x}\right)^{\alpha+1}\frac{\ud
x}{x}\qquad\textrm{is infinite}.
\]
\end{itemize}
\end{corollary}
\begin{corollary}
The upper envelope of the stable subordinator with index $\alpha\in
(0,1)$ at $0$ and at $+\infty$ is as follows:
\begin{itemize}
\item[i)] Let $h\in\mathcal{H}_{0}$, such that either $\lim_{t\to
0}t/h(t)=0$ or $\liminf_{t\to 0}t/h(t)>0$, then
\[
\p\Big(X^{(0)}_{t}>h(t), \textrm{ i.o., as } t\to 0\Big)= 0\textrm{
or } 1,
\]
according as
\[
\int_{0^+} \left(\frac{h(x)}{x}\right)^{\alpha}\frac{\ud
x}{x}\quad\textrm{is finite } \quad\textrm{or}\quad \int_{0^+}
\left(\frac{h(x)}{x}\right)^{\alpha+1}\frac{\ud
x}{x}\qquad\textrm{is infinite}.
\]
\item[ii)]Let $h\in\mathcal{H}_{\infty}$, such that either
$\lim_{t\to +\infty}t/h(t)=0$ or $\liminf_{t\to +\infty}t/h(t)>0$,
then for all $x\geq 0$
\[
\p\Big(X^{(x)}_{t}>h(t), \textrm{ i.o., as } t\to \infty\Big)=
0\textrm{ or } 1,
\]
according as
\[
\int^{+\infty}\left(\frac{h(x)}{x}\right)^\alpha\frac{\ud
x}{x}\qquad\textrm{is finite }\quad \textrm{or}\quad
\int^{+\infty}\left(\frac{h(x)}{x}\right)^{\alpha+1}\frac{\ud
x}{x}\qquad\textrm{is infinite}.
\]
\end{itemize}
\end{corollary}
\noindent\textbf{Example 2.} Let $X^{(0)}$ be a stable L\'evy
process conditioned to stay positive with no positive jumps and
index $1<\alpha\leq 2$ (See Bertoin \cite{be} for a proper
definition).\\
From Pardo \cite{Pa}, we know that $X^{(0)}$ drifts towards
$+\infty$ and that
\[
-\log \hat{F}(1/t)\sim
\frac{\alpha-1}{\alpha}\left(\frac{1}{\alpha}\right)^{1/(\alpha-1)}x^{1/\alpha-1}\quad\textrm{as}\quad
t\to +\infty.
\]
Then applying Theorems 5 and 6, and Corollary 3, we get the
following laws of the iterated logarithm.
\begin{corollary}
Let $X^{(0)}$ be a stable L\'evy process conditioned to stay
positive with no positive jumps and $\alpha>1$. Then, the first
passage time process satisfies
\[
\liminf_{t\to 0}\frac{S_{t}\big(\log |\log t|\big)^{\alpha
-1}}{t^{\alpha}}=\frac{1}{\alpha}\left(1-\frac{1}{\alpha}\right)^{\alpha-1},
\quad\textrm{ almost surely.}
\]
The same law of the iterated logarithm is satisfied for large
times.\\
The processes $X^{(x)}$, $J^{(x)}$ and $X^{(x)}-J^{(x)}$ satisfy the
following laws of the iterated logarithm
\[
\limsup_{t\to 0}\frac{X^{(0)}_{t}}{t^{1/\alpha}\big(\log |\log
t|\big)^{1-1/\alpha}}=\alpha\left(\alpha-1\right)^{-\frac{\alpha-1}{\alpha}},
\quad\textrm{ almost surely,}
\]
\[
\limsup_{t\to 0}\frac{J^{(0)}_{t}}{t^{1/\alpha}\big(\log |\log
t|\big)^{1-1/\alpha}}=\alpha\left(\alpha-1\right)^{-\frac{\alpha-1}{\alpha}},
\quad\textrm{ almost surely,}
\]
\[
\limsup_{t\to 0}\frac{X^{(0)}_{t}-J^{(0)}_{t}}{t^{1/\alpha}\big(\log
|\log
t|\big)^{1-1/\alpha}}=\alpha\left(\alpha-1\right)^{-\frac{\alpha-1}{\alpha}},
\quad\textrm{ almost surely,}
\]
and for all $x\geq 0$,
\[
\limsup_{t\to +\infty}\frac{X^{(x)}_{t}}{t^{1/\alpha}\big(\log \log
t\big)^{1-1/\alpha}}=\alpha\left(\alpha-1\right)^{-\frac{\alpha-1}{\alpha}},
\quad\textrm{ almost surely,}
\]
\[
\limsup_{t\to +\infty}\frac{J^{(x)}_{t}}{t^{1/\alpha}\big(\log \log
t\big)^{1-1/\alpha}}=\alpha\left(\alpha-1\right)^{-\frac{\alpha-1}{\alpha}},
\quad\textrm{ almost surely,}
\]
\[
\limsup_{t\to
+\infty}\frac{X^{(x)}_{t}-J^{(x)}_{t}}{t^{1/\alpha}\big(\log \log
t\big)^{1-1/\alpha}}=\alpha\left(\alpha-1\right)^{-\frac{\alpha-1}{\alpha}},
\quad\textrm{ almost surely.}
\]
\end{corollary}
\noindent{\bf Example 3.} Let $\xi$ be a L\'evy process which drifts
towards $+\infty$ and with finite exponential moments of arbitrary
positive order. Note that this condition is satisfied, for example,
when the jumps of $\xi$ are bounded from above by some fixed number,
in particular when $\xi$ is a L\'evy process with no positive jumps.
More precisely, we have
\[
\e\big(e^{\lambda
\xi_{t}}\big)=\exp\big\{t\psi(\lambda)\big\}<+\infty\qquad t,
\lambda\geq 0.
\]
From Theorem 25.3 in Sato \cite{Sa}, we know that this hypothesis is
equivalent to assume that the L\'evy measure $\Pi$ of $\xi$
satisfies
\[ \int_{[1,\infty)}e^{\lambda x}\Pi(\ud x)
<+\infty\qquad \textrm{for every  } \lambda>0.
\]
Under this condition and with the hypothesis that $\psi$, the
Laplace exponent of $\xi$, varies regularly at $+\infty$ with index
$\beta\in(1,2)$, Pardo \cite{Pa} gave the following estimates of the
tail probabilities of $I(\hat{\xi})$ and $\nu I(\hat{\xi})$,

\begin{equation*}
-\log \p\Big(\nu I(\hat{\xi})<1/x\Big)\sim-\log
\p\left(I(\hat{\xi})<1/x\right)\sim
(\beta-1)\overset{\leftharpoonup}{H}(x)\quad \textrm{ as }\quad x\to
+\infty,
\end{equation*}
where
\begin{equation*}
\overset{\leftharpoonup}{H}(x)=\inf\Big\{s > 0\textrm{ , }
\psi(s)/s> x\Big\},
\end{equation*}
is a regularly varying function with index $(\beta -1)^{-1}$.\\
Hence, the pssMp associated to $\xi$ satisify condition
($\ref{logreg}$). This allow us to obtain laws of iterated logarithm
for the first
passage time process and for the pssMp in terms of the following function.\\
Let us define the function
\[
f(t):=\frac{\log|\log t|}{\psi(\log|\log t|)}\qquad
\textrm{for}\qquad t>1,\quad t\neq e.
\]
By integration by parts, we can see that the function
$\psi(\lambda)/\lambda$ is increasing, hence it is straightforward
that the function $tf(t)$ is also increasing in a neighbourhood of
$\infty$.
\begin{corollary} If $\psi$ is regularly varying at $+\infty$ with index $\beta\in
(1,2)$, then
\[
\liminf_{x\to 0}\frac{S_x}{xf(x)}=(\beta-1)^{\beta -1} \qquad
\textrm{almost surely}
\]
and,
\begin{equation*}
\liminf_{x\to +\infty}\frac{S_x}{xf(x)}=(\beta-1)^{\beta -1} \qquad
\textrm{almost surely.}
\end{equation*}
\end{corollary}
Let us define
\[
g(t):=\frac{\psi(\log |\log t|)}{\log |\log t|}\qquad{\rm for}\quad
t>1, \quad t\neq e.
\]
\begin{corollary}If $\psi$ is regularly varying at $+\infty$ with index $\beta\in
(1,2)$, then
\begin{equation*}
\limsup_{t\to 0}\frac{X_{t}}{tg(t)}=(\beta-1)^{-(\beta -1)} \qquad
\textrm{almost surely}
\end{equation*}
and for all $x\geq 0$,
\begin{equation*}
\limsup_{t\to +\infty}\frac{X^{(x)}_t}{tg(t)}=(\beta-1)^{-(\beta
-1)} \qquad \textrm{almost surely}.
\end{equation*}
Moreover, if the processes $X^{(x)}$ has no positive jumps,
$J^{(x)}$ and $X^{(x)}-J^{(x)}$ satisfy the following laws of the
iterated logarithm:
\[
\limsup_{t\to 0}\frac{J^{(0)}_{t}}{tg(t)}=(\beta-1)^{-(\beta-1)},
\quad\textrm{ almost surely,}
\]
\[
\limsup_{t\to
0}\frac{X^{(0)}_{t}-J^{(0)}_{t}}{tg(t)}=(\beta-1)^{-(\beta-1)},
\quad\textrm{ almost surely,}
\]
and for all $x\geq 0$,
\[
\limsup_{t\to
+\infty}\frac{J^{(x)}_{t}}{tg(t)}=(\beta-1)^{-(\beta-1)},
\quad\textrm{ almost surely,}
\]
\[
\limsup_{t\to
+\infty}\frac{X^{(x)}_{t}-J^{(x)}_{t}}{tg(t)}=(\beta-1)^{-(\beta-1)},
\quad\textrm{ almost surely.}
\]
\end{corollary}

\noindent\textbf{Example 4.}  Sato \cite{Sa1} ( see also Sato
\cite{Sa}) studied some interesting properties of positive
increasing self-similar processes with independent increments. In
particular, the author showed that if $Y=(Y_{t}, t\geq 0)$ is a
process with such characteristics and starting from $0$, we can
represent its Laplace transform by
\[
\e\Big[\exp\big\{-\lambda
Y_1\big\}\Big]=\exp\big\{-\kappa(\lambda)\big\}\qquad\textrm{for
}\lambda>0,
\]
where
\[
\kappa(\lambda)=c\lambda+\int_0^{+\infty}\big(1-e^{-\lambda
x}\big)\frac{l(x)}{x}\ud x,
\]
$c\geq 0$ and $l(x)$ is a nonnegative decreasing function on $(0,
+\infty)$ with
\[\int_0^{+\infty}\frac{l(x)}{1+x}\ud x
<+\infty.
\]
From its definition, it is clear that the Laplace exponent $\kappa$
is an increasing continuous function and more precisely a concave
function.\\
Under the assumption that $\kappa$ varies regularly  at
$+\infty$ with index $\alpha\in(0,1)$, we will have the
following sharp estimate for the distribution of $Y_{1}$.\\
Let us define the function
\[
\rho(t)=\frac{t\log |\log
t|}{\overset{\leftharpoonup}\kappa(\log|\log t|)},\quad \textrm{ for
}\quad t\neq e, \quad t>1,
\]
where $\overset{\leftharpoonup}\kappa$ is the inverse function of
$\kappa$.
\begin{proposition}
Let $(Y_{t}, t\geq 0)$ be a  positive increasing self-similar
processes with independent increments and suppose that $\kappa$, its
Laplace exponent, varies regularly at $+\infty$ with index
$\alpha\in(0,1)$. Then for every $c>0$, we have
\[
-\log\p\left(Y_{1}\leq
\frac{c\rho(t)}{t}\right)\sim\big(\alpha/c\big)^{\frac{\alpha}{1-\alpha}}(1-\alpha)\log|\log
t|\quad\textrm{ as}\quad t\to 0 \quad(t\to \infty).
\]
\end{proposition}
\textit{Proof:} From de Brujin's Tauberian Theorem (see for instance
Theorem 4.12.9 in \cite{Bing}), we have that if $\kappa$ varies
regularly at $+\infty$ with index $\alpha\in(0,1)$ then
\[
-\log \p\big(Y_{1}\leq x\big)\sim
\frac{\alpha^{\frac{\alpha}{1-\alpha}}(1-\alpha)}{\overset{\leftharpoonup}\omega(1/x)},
\quad\textrm{ for }\quad x\to 0,
\]
where $\overset{\leftharpoonup}\omega$ is the asymptotic inverse of
$\omega$, a regularly varying function at $+\infty$ with index
$(\alpha -1)/\alpha$ and that satisfies
\begin{equation}\label{equprop5}
\frac{\lambda}{\kappa(\lambda)}\sim
\omega\left(\frac{1}{\kappa(\lambda)}\right) \quad \textrm{for }
\quad \lambda \to +\infty.
\end{equation}
Hence, taking $x=c\rho(t)/t$ and
$\lambda=\overset{\leftharpoonup}\kappa(\log|\log t|)$, and doing
some calculations we get the desired result.\QED This estimate allow
us to get the following law of the iterated logarithm.
\begin{corollary} Let $(Y_t, t\geq 0)$ be a positive increasing self-similar
processes with independent increments and suppose that $\kappa$, its
Laplace exponent, satisfies the conditions of the previous
Proposition. Then, we have
\[
\liminf_{t\to
0}\frac{Y_t}{\rho(t)}=\alpha(1-\alpha)^{(1-\alpha)/\alpha}, \qquad
\textrm{almost surely.}
\]
The same law of iterated logarithm is satisfied for large times.
\end{corollary}
\textit{Proof}: The proof of this Corollary follows from the
integral test found by Watanabe \cite{wa} and applying the above
estimate of the distribution of $Y_1$.\QED

 Now, we denote by $\kappa_{1}$ and $\kappa_{2}$ the
Laplace exponents of the last and first passage time processes,
respectively.  Since $S_1\leq U_1$, it is clear that
$\kappa_{2}(\lambda)\leq \kappa_{1}(\lambda)$
for all $\lambda \geq 0$.\\
Let us suppose that $\kappa_1$ and $\kappa_2$ are regularly varying
at $+\infty$ with index $\alpha_{1}$ and $\alpha_{2}$ respectively,
such that $0< \alpha_2\leq \alpha_1 <1$. By Proposition 9 and
Theorem 8, we can deduce that  $\kappa_1$ and $\kappa_2$ are
asymptotically equivalents and that $\alpha_{1}=\alpha_{2}$. Then by
the above Corollary, we have that for
\[
h_1(t):=\frac{t\log |\log
t|}{\overset{\leftharpoonup}\kappa_{1}(\log|\log t|)}\quad
\textrm{and}\quad h_2(t):=\frac{t\log |\log
t|}{\overset{\leftharpoonup}\kappa_{2}(\log|\log t|)},\quad \textrm{
for }\quad  t\neq e, \quad t>1,
\]
where $\overset{\leftharpoonup}\kappa_{1}$ and
$\overset{\leftharpoonup}\kappa_{2}$ are the inverse of $\kappa_1$
and $\kappa_2$, respectively; the processes $U$ and $S$ satisfy
\[
\liminf_{t\to
0}\frac{U_t}{h_{1}(t)}=\alpha_1(1-\alpha_1)^{(1-\alpha_{1})/\alpha_{1}}
\qquad \textrm{almost surely,}
\]
and
\[
\liminf_{t\to
0}\frac{S_t}{h_{1}(t)}=\alpha_1(1-\alpha_1)^{(1-\alpha_{1})/\alpha_{1}}
\qquad \textrm{almost surely.}
\]
Note that we can replace $h_{1}$ by $h_{2}$ and that we also have
the same laws of the iterated logarithm for
large times.\\
By the sharp estimation in Proposition 10 of the tail probability of
$S_{1}$, we deduce the following law of the
iterated logarithm.\\
Let us define
\[
f_{2}(t):=\frac{t\overset{\leftharpoonup}\kappa_{2}(\log|\log
t|)}{\log|\log t|},\quad \textrm{ for } \quad t\neq e,\quad t>1.
\]
\begin{corollary} Let $\kappa_{2}$ be the Laplace exponent of $S_{1}$.
If $\kappa_{2}$ is regularly varying at $+\infty$ with index
$\alpha_{2}\in(0,1)$, then
\[
\limsup_{t\to
0}\frac{X^{(0)}_t}{f_{2}(t)}=\alpha_{2}^{-1}(1-\alpha_{2})^{-(1-\alpha_{2})/\alpha_{2}}\qquad
\textrm{ almost surely,}
\]
and for any $x\geq 0$,
\[
\limsup_{t\to
+\infty}\frac{X^{(x)}_t}{f_{2}(t)}=\alpha_{2}^{-1}(1-\alpha_{2})^{-(1-\alpha_{2})/\alpha_{2}}\qquad
\textrm{ almost surely.}
\]
\end{corollary}
On the other hand, from Theorem 8 we get the following Corollary.
\begin{corollary}
Let $\kappa_{2}$ be the Laplace exponent of $S_{1}$ and
$\overset{\leftharpoonup}\kappa_{2}$ its inverse. If $\kappa_{2}$ is
regularly varying at $+\infty$ with index $\alpha_{2}\in(0,1)$, then
\begin{itemize}
\item[i)]\[ \limsup_{t\to
0}\frac{J^{(0)}_t}{f_{2}(t)}=\alpha_{2}^{-1}(1-\alpha_{2})^{-(1-\alpha_{2})/\alpha_{2}}\qquad
\textrm{ almost surely,}
\]
and for any $x\geq 0$,
\[
\limsup_{t\to
+\infty}\frac{J^{(x)}_t}{f_{2}(t)}=\alpha_{2}^{-1}(1-\alpha_{2})^{-(1-\alpha_{2})/\alpha_{2}}\qquad
\textrm{ almost surely.}
\]
\item[ii)]
\[
\limsup_{t\to
0}\frac{X^{(0)}_t-J^{(0)}_t}{f_{2}(t)}=\alpha_{2}^{-1}(1-\alpha_{2})^{-(1-\alpha_{2})/\alpha_{2}}\qquad
\textrm{ almost surely,}
\]
and for any $x\geq 0$,
\[
\limsup_{t\to
+\infty}\frac{X^{(x)}_{t}-J^{(x)}_t}{f_{2}(t)}=\alpha_{2}^{-1}(1-\alpha_{2})^{-(1-\alpha_{2})/\alpha_{2}}\qquad
\textrm{ almost surely.}
\]
\end{itemize}
\end{corollary}

Now, we will apply  these results to the case of transient Bessel
process. Here, we employ the usual Bessel functions $I_{a}$ and
$K_{a}$, as in Kent \cite{Ke} and Jeanblanc, Pitman and Yor
\cite{JPY}. It is well known that
\[
\e\left(\exp\Big\{-\lambda
S_{1}\Big\}\right)=\lambda^{a/2}\frac{1}{2^{a/2}\Gamma(a+1)I_{a}(\sqrt{2\lambda})},\quad
\lambda>0,
\]
and
\[
\e\left(\exp\Big\{-\lambda
U_{1}\Big\}\right)=\frac{\lambda^{a/2}}{2^{
a/2-1}\Gamma(a)}K_{a}(\sqrt{2\lambda}),\quad \lambda>0,
\]
where $\Gamma$ is the gamma function (see for instance Jeanblanc, Pitman and Yor \cite{JPY}).\\
Now, we define for $\lambda>0$
\[
\phi_{1}(\lambda)=\log(2^{a/2-1}\Gamma(a))-\log
K_{a}(\sqrt{2\lambda}) -\log \lambda^{a/2},
\]
\[
\phi_{2}(\lambda)=\log I_{a}(\sqrt{2\lambda})
+\log(2^{a/2}\Gamma(a+1)) -\log \lambda^{a/2}.
\]
Since, we have the following asymptotic behaviour
\[
I_{a}(x)\sim (2\pi x)^{-1/2}e^{x}\quad \textrm{and}\quad
K_{a}(x)\sim\left(\frac{\pi}{2x}\right)^{1/2}e^{-x}\quad \textrm{
when }x\to +\infty,
\]
(see Kent \cite{Ke} for instance), we deduce that $\phi_1$ and
$\phi_{2}$ are regularly varying at $+\infty$ with index $1/2$. From
Proposition 9 and Theorem 8, we deduce that they are asymptotically
equivalent.\\
From the above results, we have that
\[
\liminf_{t\to 0}\frac{U_{t}}{h_{1}(t)}=1/4, \qquad \liminf_{t\to
0}\frac{S_{t}}{h_{2}(t)}=1/4 \quad\textrm{almost surely,}
\]
\[
\limsup_{t\to 0}\frac{X^{(0)}_{t}}{f_{2}(t)}=4, \qquad \limsup_{t\to
0}\frac{J^{(0)}_{t}}{f_{1}(t)}=4 \quad\textrm{almost surely,}
\]
and
\[
\limsup_{t\to 0}\frac{J^{(0)}_{t}}{f_{2}(t)}=4, \qquad \limsup_{t\to
0}\frac{X^{(0)}_{t}-J^{(0)}_{t}}{f_{1}(t)}=4 \quad\textrm{almost
surely,}
\]
where
\[
h_{1}(t)=\frac{t\log |\log
t|}{\overset{\leftharpoonup}\kappa_{1}(\log|\log t|)},\quad
h_{2}(t)=\frac{t\log |\log
t|}{\overset{\leftharpoonup}\kappa_{2}(\log|\log t|)},\quad
f_{1}(t)=\frac{t^{2}}{h_{1}(t)}, \quad
f_{2}(t)=\frac{t^{2}}{h_{2}(t)}
\]
and, $\overset{\leftharpoonup}\kappa_{1}$ and
$\overset{\leftharpoonup}\kappa_{2}$ are the inverse functions of
$\kappa_{1}$ and $\kappa_{2}$, respectively.\\
Similarly, we have all these laws of the iterated logarithm for
large times.\\
\noindent\textbf{Example 5.} Let $\xi=N$ be a standard Poisson
process. From Proposition 3 in Bertoin and Yor \cite{BeY2} and
Example 1 in Pardo \cite{Pa}, we know
\[
-\log\p\big(I(\hat{\xi})<t\big)\sim -\log\p\big(\nu
I(\hat{\xi})<t\big)\sim \frac{1}{2}(\log 1/t)^2,\qquad\textrm{as }
t\to 0.
\]
Hence, we obtain the following laws of the iterated logarithm. Let
us define
\[
m(t):=t\exp\Big\{-\sqrt{2\log|\log t|}\Big\}.
\]
\begin{corollary}
Let $N$ be a Poisson process, then the pssMp $X^{(x)}$ associated to
$N$ by the Lamperti representation satisfies the following law of
the iterated logarithm,
\[
\limsup_{t\to 0}\frac{X^{(0)}_t m(t)}{t^2}=1\qquad \textrm{almost
surely.}
\]
For $x\geq 0$
\[
\limsup_{t\to +\infty}\frac{X^{(x)}_t m(t)}{t^2}=1\qquad
\textrm{almost surely.}
\]
The first passage time process $S$ associated to $X^{(0)}$ satisfies
the following law of the iterated logarithm,
\[
\liminf_{x\to 0}\frac{S_x}{m(x)}=1,\quad\textrm{ and }\quad
\liminf_{x\to +\infty}\frac{S_x}{m(x)}=1, \quad{ almost surely.}
\]
\end{corollary}

\noindent {\bf Acknowledgements.} I would like to thank Lo\"\i c
Chaumont for guiding me through the development of this work, and
for all his helpful advice.

\end{document}